\documentclass[twoside,twocolumn]{article}
\usepackage{acta}
\usepackage{graphicx}
\usepackage{dcolumn}
\usepackage{amsmath,amsfonts,amssymb}
\usepackage{bm}
\newcommand{\D}{\partial}
\newcommand{\DD}{\delta}
\newcommand{\F}{\mathcal{F}}
\newcommand{\frakl}{\mathfrak{l}}
\newcommand{\I}{\mathbf{I}}
\begin{document}
\Eingang{27}{10}{2012} \Annahme{6}{12}{2012} \Sachgebiet{1}
\Band{98} \Jahr{2012} \Heft{} \Ersteseite{1} \Letzteseite{}

\AuthorsI{S. Bilbao$^{1)}$, A. Torin$^{1)}$, V. Chatziioannou$^{2)}$}
\AddressI{$^{1)}$ Acoustics and Audio Group, King's Buildings, University of Edinburgh, Edinburgh, United Kingdom.\\
\hspace*{8pt}sbilbao@staffmail.ed.ac.uk\\
$^{2)}$ Institute of Music Acoustics, University of Music and Performing Arts, Vienna, Austria}


\Englishtitle{Numerical Modeling of Collisions in Musical Instruments} \Germanorfrenchtitle{}

\Kolumnentitel{Bilbao, Torin, Chatziioannou: Numerical Modeling of Collisions}

\Englishabstract{Collisions play an important role in many aspects of the physics of musical instruments. The striking action of a hammer or mallet in keyboard and percussion instruments is perhaps the most important example, but others include reed-beating effects in wind instruments, the string/neck interaction in fretted instruments such as the guitar as well as in the sitar and the wire/membrane interaction in the snare drum. From a simulation perspective, whether the eventual goal is the validation of musical instrument models or sound synthesis, such highly nonlinear problems pose various difficulties, not the least of which is the risk of numerical instability. In this article, a novel finite difference time domain simulation framework for such collision problems is developed, where numerical stability follows from strict numerical energy conservation or dissipation, and where a a power law formulation for collisions is employed, as a potential function within a Hamiltonian formulation. The power law serves both as a model of deformable collision, and as a mathematical penalty under perfectly rigid, non-deformable collision. This formulation solves a major problem underlying previous work, where a Hamiltonian framework was not employed for collisions, and thus stability was not ensured. Various numerical examples, illustrating the unifying features of such methods across a wide variety of systems in musical acoustics are presented, including numerical stability and energy conservation/dissipation, bounds on spurious penetration in the case of rigid collisions, as well as various aspects of musical instrument physics.}
\Germanorfrenchabstract{}

\ScientificPaper

\section{Introduction}

Various mechanisms of sound production in musical instruments rely on collisions; the obvious examples are the interaction of a striking object, such as a hammer with a string, or a drum stick or mallet with a percussive instrument, but more subtle examples include reed-beating effects in wind instruments, string/neck interactions along the fretboard of a guitar, and also in the sitar and tambura, and most dramatically, the wire/membrane interaction in a snare drum. In all cases, the collision interaction is necessarily strongly nonlinear, and simulation design becomes a challenging problem. 

Such collisions may be grouped into two types. In the first, one of the two objects involved in the collision is modelled as lumped---i.e., it is characterized by a single position/velocity pair. Such is the case for most models of the hammer-string \cite{Boutillon88, Chaigne}, mallet-bar \cite{Chaigne97} and mallet-membrane interactions \cite{Rhaouti}, even if the colliding object occupies a finite interaction region. The perceptual effects of such a nonlinear interaction are major, particularly with regard to the spectral content of the resulting sound. In the second, both objects must be considered to be fully distributed, and the region of contact will vary in a non-trivial manner. Examples include the snare drum \cite{Bilbaosnare}, sitar \cite{vyasarayani09} and string fret interactions \cite{evangelista11} mentioned above. 

Collision modeling is a branch of the large area of contact mechanics with applications across a variety of disciplines, and especially robotics \cite{Marhefka99} and computer graphics \cite{baraff94}---see Wriggers et al.\cite{wriggers08} for a recent review. In musical acoustics and applications in sound synthesis, various techniques have been employed, including digital waveguides \cite{krishnaswamy03, evangelista11}, modal methods \cite{vyasarayani09} and time stepping methods in the lumped setting \cite{avanzini01, papetti11, chatziioannou13}, in modeling lumped/distributed collisions \cite{Chabassier, Chaigne, Bilbaonewbook}, the interaction of a distributed object with a rigid barrier \cite{krishnaswamy03, Siddiq12, kartofelev13}, and in the collision of deformable objects for synthesis applications\cite{bensoam10}. Finite element methods are often used to model collisions of complex deformable bodies in mainstream applications\cite{hughes76}; in the present case of musical acoustics and sound synthesis, where geometries are often simple, finite difference time domain methods \cite{Strikwerda} are an efficient alternative, and will be employed here.

In the setting of musical acoustics, previous work on particular cases of collisions using finite difference methods, such as, for example, the snare membrane interaction \cite{Bilbaosnare} and the reed/lay interaction \cite{Bilbaonewbook} has generally been dealt with using ad hoc techniques---that is to say, there is not any attempt at proving numerical stability under such strongly nonlinear conditions. As in the case of other types of inherent distributed nonlinearities in musical instruments (as, e.g., in strings \cite{Bankjasa, Bilbaojasa05, Chabassier} or curved shells \cite{Thomas04}), an approach based on energy principles is of great utility; energetic methods are widely used in elastodynamics, whether in a variational setting \cite{vouga10} or when strict energy conservation or dissipation is enforced \cite{simo92}. The general aim of this paper, then, is to incorporate the nonlinear collision mechanism into the energy conservation framework, and use such a formulation in order to arrive at sufficient stability conditions for such methods, and show its application to a wide variety of systems, spanning the range of musical instruments. The mechanism by which this is accomplished is through the introduction of a potential energy term corresponding to the collision interaction.

If the collision is accompanied by some deformation of the colliding object (as for, e.g., the piano hammer), then the additional potential has the interpretation of energy stored in the deformed object. In the case of rigid collisions, however, the potential must be interpreted as a {\em penalty}---such a potential results in a strong opposing force, penalizing penetration, but some limited interpenetration of the colliding bodies is permitted. Such penalty methods \cite{Harmonthesis} can be viewed in contrast with methods based on hard non-penetrative constraints \cite{fetecau03}. As long as contact velocities are low (which is generally the case in musical acoustics applications), the latter approach is justifiable. Indeed, a second benefit of such an energy-based formulation, beyond proving stability, is a means of determining bounds on such spurious penetration---as will be seen, the amount of penetration can be made as small as desired (and particular, small enough to be negligible in applications in acoustics).  Numerical existence and uniqueness results in the simple case of a lumped collision with a rigid barrier have been presented recently \cite{chatziioannou13} and will be extended to fully distributed systems here. 

In Section \ref{lumpedsec}, the collision of a mass with a rigid barrier is used as a test case in the construction of stable energy conserving methods, accompanied by an extension to the case of a lossy collision. Collisions of a lumped object with a distributed system are considered next. To this end, some background material on discrete representations of distributed systems is presented in Section \ref{backgdsec}.  The case of the hammer string interaction is treated first in Section \ref{hammersec}, followed by reed beating effects, in a simple model of a wind instrument in Section \ref{reedsec}. A fully distributed collision is considered next, in the case of a string in contact with a rigid barrier in Section \ref{barriersec}. In 2D, the simplest system of interest is the mallet/membrane interaction, described in Section \ref{malletsec}. Finally, the case of a wire in contact with a membrane, modeling the snare interaction is described in Section \ref{snaresec}. As a nontrivial example, in this last section, a complete 3D model of a snare drum, including the cavity, membranes, acoustic field and snares is outlined. As the schemes presented here all require the solution of nonlinear equations, some existence/uniqueness results are provided in Appendix \ref{newtonsec}, along with some comments on the use of iterative methods.

\section{Prelude: Collision of a Mass with a Rigid Barrier}
\label{lumpedsec}
A starting point in many numerical studies of collisions \cite{wriggers08, papetti11, chatziioannou13} is the case of a lumped deformable object approaching a rigid barrier from below, under a nonlinear interaction force. The system may be written as
\begin{equation}
\label{lumpeddef}
M\frac{d^2 u}{dt^2} = -f \qquad f = \frac{d\Phi}{du} = \frac{d\Phi}{dt}/\frac{du}{dt} 
\end{equation}
Here, $u=u(t)$ is the position of the object at time $t$, and $M$ is its mass. $f=f(t)$ is the interaction force, written here in terms of a potential $\Phi(u)\geq 0$. For collisions, with a barrier at $u=0$, $f$ is zero for $u\leq 0$; when $u>0$, the object undergoes compression. A power law is a useful general choice of the potential:
\begin{equation}
\label{phidef}
\hspace{-0.2in}\Phi = \Phi_{K,\alpha} = \frac{K[u]_{+}^{\alpha+1}}{\alpha+1}\geq 0 \quad\rightarrow\quad f = K[u]_{+}^{\alpha}
\end{equation}
It depends on a stiffness parameter $K\geq 0$ and an exponent $\alpha> 1$; in this article, the notation $[\cdot]_{+}$ indicates the positive part of, i.e., $[u]_{+} = \frac{1}{2}\left(u+|u|\right)$.

Such power law nonlinearities are common in models of collisions in many settings, including not just the present case of the lumped collision, but also in models of the hammer string interaction \cite{Chaigne}, where $\alpha$ and $K$ are empirically determined (though for simpler systems, such as the contact of two spheres, Hertzian models allow a direct calculation \cite{Horvay80}). Such models may be viewed as permitting a certain deformation of the colliding object when in contact with the barrier. If the colliding object is perfectly rigid, then such a model is obviously unphysical, in that some penetration is permitted---in this case, the potential $\Phi$ may be interpreted as a penalty; in numerical treatments, $\alpha=1$ is often chosen \cite{Harmonthesis}. Most of what follows here does not depend on a particular choice of the potential---but existence and uniqueness in implementation do depend on the form of $\Phi$, and the power law has various advantages which will be highlighted.

\subsection{Energy Balance}

The collision model presented above is energy conserving \cite{papetti11}. Multiplying \eqref{lumpeddef} by $du/dt$ leads to 
\begin{equation}
M\frac{du}{dt}\frac{d^2 u}{dt^2} = -\frac{du}{dt}f = -\frac{du}{dt}\frac{d\Phi}{du} = -\frac{d\Phi}{dt}
\end{equation}
and thus to the energy balance, in terms of total energy ${\mathcal H}$:
\begin{equation}
\label{egybal}
\frac{d{\mathcal H}}{dt} = 0 \qquad {\mathcal H}(t) = \frac{M}{2}\left(\frac{du}{dt}\right)^2+\Phi
\end{equation}
Thus ${\mathcal H}(t)$, corresponding to the total energy of the object at time $t$, is conserved and non-negative:
\begin{equation}
{\mathcal H}(t) = {\mathcal H}\left(0\right)\geq 0
\end{equation}

For the particular choice of the power law potential function $\Phi = \Phi_{K,\alpha}$ from \eqref{phidef}, the conservation law implies bounds on both $u(t)$ and $du/dt$ at any time $t$, in terms of the initial energy ${\mathcal H}(0)$:
\begin{equation}
\hspace{-0.2in}\Big|\frac{du}{dt}\Big|\leq \sqrt{\frac{2{\mathcal H}(0)}{M}}\qquad u(t)\leq \left(\frac{\left(\alpha+1\right){\mathcal H}(0)}{K}\right)^{\frac{1}{\alpha+1}}
\end{equation}
The first bound holds for any non-negative potential $\Phi$, and the second employs monotonicity of $\Phi_{K,\alpha}$ in $u\geq 0$. 

\subsection{Time Series and Difference Operators}
\label{timeopsec}
Moving to a discrete time simulation setting, let $u = u^{n}$ represent an approximation to $u(t)$ at $t=nk$, for integer $n$, and for a given time step $k$.

Unit shifts $e_{t+}$ and $e_{t-}$ are defined as
\begin{equation}
e_{t+}u^{n} = u^{n+1}\qquad e_{t-}u^{n} = u^{n-1}
\end{equation}
Forward, backward and centered difference approximations to a first time derivative may thus be defined as
\begin{equation}
\label{dtdef}
\hspace{-0.3in}\delta_{t+} = \frac{e_{t+}-1}{k}\quad \delta_{t-} = \frac{1-e_{t-}}{k}\quad \delta_{t\cdot} = \frac{e_{t+}-e_{t-}}{2k}
\end{equation}
and an approximation to a second time derivative as
\begin{equation}
\label{dttdef}
\delta_{tt} = \frac{e_{t+}-2+e_{t-}}{k^2}
\end{equation}

Various averaging operators may be defined as
\begin{equation}
\label{mudef}
\hspace{-0.3in}\mu_{t+} = \frac{e_{t+}+1}{2}\quad \mu_{t-} = \frac{1+e_{t-}}{2}\quad \mu_{t\cdot} = \frac{e_{t+}+e_{t-}}{2}
\end{equation}

\subsection{Finite Difference Scheme}

 An approximation to \eqref{lumpeddef}, at time step $n$ is then
\begin{subequations}
\label{lumpedtotalfd}
\begin{equation}
\label{lumpedfddef}
M\delta_{tt}u^{n} = -f^{n}
\end{equation}
where $\delta_{tt}$ is a second difference operator, as defined in \eqref{dttdef}, and where $f^{n}$ is a time series defined by
\begin{equation}
\label{fddef}
f^{n} = \frac{\delta_{t-}\Phi^{n+\frac{1}{2}}}{\delta_{t\cdot}u^{n}}\qquad \Phi^{n+\frac{1}{2}} = \mu_{t+}\Phi\left(u^{n}\right)
\end{equation}
\end{subequations}
in terms of a discrete potential $\Phi^{n+\frac{1}{2}}$, itself defined in terms of the potential $\Phi$ of the model problem. Note here that the definition of $f^{n}$ here mirrors that of the continuous time case, from \eqref{lumpeddef}; simpler, fully explicit forms are available, but can lead to stability problems---see the end of Section \ref{hammersimsec} for an example of numerical instability, in the case of a hammer in contact with a string. The backwards and centered difference operators $\delta_{t-}$ and $\delta_{t\cdot}$ employed in \eqref{fddef} are defined in \eqref{dtdef}, and the averaging operator $\mu_{t+}$ is defined in \eqref{mudef}. The time series $\Phi^{n+\frac{1}{2}}$ is interleaved with respect to the time series $u^{n}$ itself; thus scheme \eqref{lumpedtotalfd} is centered about time step $n$ and, if stable, is thus second order accurate\cite{Strikwerda}. 

Multiplying \eqref{lumpedfddef} by $\delta_{t\cdot}u^{n}$, and employing \eqref{fddef} gives
\begin{equation}
M\delta_{t\cdot}u^{n}\delta_{tt}u^{n} = -\delta_{t\cdot}u^{n}f^{n} = -\delta_{t-}\Phi^{n+\frac{1}{2}}
\end{equation}
Employing the identity
\begin{equation}
\delta_{t\cdot}u^{n}\delta_{tt}u^{n} = \frac{1}{2}\delta_{t-}\left(\delta_{t+}u^{n}\right)^2
\end{equation}
leads to the discrete energy balance
\begin{equation}
\label{egybalfd}
\delta_{t-}{\mathfrak h}^{n+\frac{1}{2}} = 0\quad\longrightarrow\quad {\mathfrak h}^{n+\frac{1}{2}} = {\rm const.}\geq 0
\end{equation}
where the discrete energy of the system is defined as
\begin{equation}
{\mathfrak h}^{n+\frac{1}{2}} = \frac{M}{2}\left(\delta_{t+}u^{n}\right)^2+\Phi^{n+\frac{1}{2}} \geq 0
\end{equation}
Bounds on the solution size follow as in the continuous case, under the choice of a power law nonlinearity $\Phi_{K,\alpha}$: 
\begin{equation}
\label{lumpedbound}
\hspace{-0.2in}|\delta_{t+}u^{n}|\leq \sqrt{\frac{2{\mathfrak h}}{M}}\qquad u^{n}\leq  \left(\frac{2\left(\alpha+1\right){\mathfrak h}}{K}\right)^{\frac{1}{\alpha+1}}
\end{equation}
and the scheme is thus unconditionally stable. Notice in particular that the amount of penetration may be controlled through the choice of $K$ in this case. 

The scheme \eqref{lumpedtotalfd} requires the solution of a nonlinear equation at each time step,
\begin{equation}
\label{Gdef}
G(r) = r+\frac{m}{r}\left(\Phi(r+a)-\Phi(a)\right)+b = 0
\end{equation}
in terms of an unknown $r$, defined in terms of $u^{n}$ as
\begin{equation}
r = u^{n+1}-u^{n-1}
\end{equation}
where $m$, $a$ and $b$ are given by
\begin{equation}
\hspace{-0.2in}m = k^2 / M\qquad a = u^{n-1}\qquad b = -2u^{n}+2u^{n-1}
\end{equation}
At each time step, once $r$ is determined, the next displacement $u^{n+1}$ may be determined from $r$ and $u^{n-1}$. 

The function $G(r)$, as defined in \eqref{Gdef}, as well as certain generalized forms, play a central role in all the algorithms to be described subsequently here in the distributed case. It has recently been shown\cite{chatziioannou13} that $G(r)=0$ possesses a unique solution for one-sided potentials of power law form. An iterative method such as the Newton-Raphson algorithm may be employed, but global convergence is not obvious. See Appendix \ref{newtonsec} for more discussion.  

\subsection{Nonlinear Losses}
\label{nllosssec}
In some more refined collision models\cite{papetti11}, a nonlinear damping term is included. From the model of Hunt and Crossley\cite{Hunt} the force $f$ in \eqref{lumpeddef} may be augmented to
\begin{equation}
\label{lossyforcedef}
f = \frac{d\Phi}{du}+\frac{du}{dt}\Xi\left(u\right)
\end{equation}
where $\Phi$ is defined as in \eqref{phidef}, and for some function $\Xi\left(u\right)\geq 0$; in particular, the choice of 
\begin{equation}
\Xi_{K,\alpha}\left(u\right) = K\beta [u]_{+}^{\alpha}
\end{equation}
for some $\beta\geq 0$ has been employed in previous studies of collisions. (Other models used in musical acoustics, particularly in the case of hammer felt are similar\cite{Stulov04}.) In this case, the energy balance \eqref{egybal} can be generalized to
\begin{equation}
\label{egyballoss}
\frac{d{\mathcal H}}{dt} = -{\mathcal Q} \quad{\rm ,}\quad {\mathcal Q} = \left(\frac{du}{dt}\right)^2\Xi(u)\geq 0
\end{equation}
and thus energy is monotonically decreasing, and the bounds on the velocity and displacement hold as before.

Scheme \eqref{lumpedtotalfd} may be generalized using
\begin{equation}
\label{fdlossdef}
f^{n} = \frac{\delta_{t-}\Phi^{n+\frac{1}{2}}}{\delta_{t\cdot}u^{n}}+\delta_{t\cdot}u^{n}\Xi^{n}
\end{equation}
where $\Xi^{n} = \Xi(u^{n})$. Now the energy balance \eqref{egybalfd} may be generalized to
\begin{equation}
\delta_{t-}{\mathfrak h}^{n+\frac{1}{2}} = -{\mathfrak q}^{n}\quad{\rm ,}\quad {\mathfrak q}^{n}= \left(\delta_{t\cdot}u^{n}\right)^2\Xi^{n}\geq 0
\end{equation}
and thus energy is monotonically decreasing, and solution bounds hold as before. The nonlinear equation to be solved at each time step must be generalized to
\begin{equation}
\label{Gdef2}
\hspace{-0.3in}G(r) = \left(1+c\right)r+\frac{m}{r}\left(\Phi(r+a)-\Phi(a)\right)+b = 0
\end{equation}
where $r$, the unknown, as well as the constants $m$, $a$ and $b$ are as before, and where $c = k\Xi(u^{n})/2M$.

\subsection{Simulations}
\label{lumpedsimsec}
\begin{figure}[htbp]
\includegraphics[width=\columnwidth]{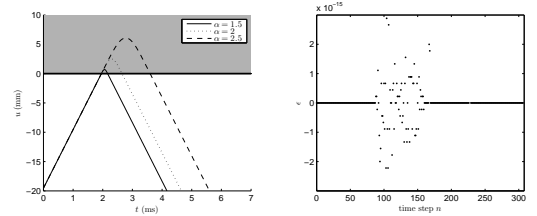}
\caption{Left: trajectory of a mass $M=10$ g, colliding with a rigid barrier, under different values of the exponent $\alpha$, as indicated. Here the mass has initial velocity 10 m/s, and $K=10^{8}$. Scheme \eqref{lumpedtotalfd} is used, with a time step $k=1/44100$ s. Right: normalized energy variation $\epsilon^{n} = \left({\mathfrak h}^{n+\frac{1}{2}}-{\mathfrak h}^{1/2}\right)/{\mathfrak h}^{1/2}$, when $\alpha = 2.5$.}
\label{lumpedfig}
\end{figure}

Consider first the collision of a mass with a rigid barrier under different choices of the exponent $\alpha$, using scheme \eqref{lumpedtotalfd}, as illustrated in Figure \ref{lumpedfig}. The trajectories of the mass are illustrated  at left; the obvious effect of the exponent is on the duration of the contact time, but as the scheme is exactly lossless, the initial and exiting speeds of the mass are identical to machine accuracy. As a further indication of the exact energy conservation property of the scheme, the normalized energy variation is plotted against time step in Figure \ref{lumpedfig}, at right; single bit variation of the energy (at double precision floating point) is evident as quantization. 

\begin{figure}[htbp]
\includegraphics[width=\columnwidth]{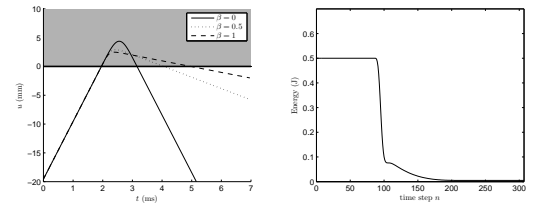}
\caption{Left: trajectory of a mass colliding with a rigid barrier, under different values of the loss coefficient $\beta$, as indicated, with $\alpha = 2.3$. All other conditions are as in the caption to Figure \ref{lumpedfig}. Right: energy ${\mathfrak h}^{n+\frac{1}{2}}$, when $\beta = 1$.}
\label{lumpedlossfig}
\end{figure}

As a further example, the behaviour of the same system, under the lossy collision model given in Section \ref{nllosssec} is plotted in Figure \ref{lumpedlossfig}, under difference choices of the loss parameter $\beta$. Energy is monotonically decreasing, illustrating a complex non-exponential loss characteristic during the collision itself. 

In the case of a completely rigid collision, some spurious penetration is inevitable. Considering again the system under the choice of parameters given in Figure \ref{lumpedfig}. Using values of $K=10^{16}$ and $\alpha = 1.2$, approximating such a  rigid collision, the penetration is under $8\times 10^{-8}$ m, which is negligible in any acoustics application.

\section{Background: Distributed Systems and Grid Functions}
\label{backgdsec}
In this section, some background material on systems and grid functions suitable for use in finite difference time domain schemes in 1D and 2D is presented, with a focus on energy techniques. 

\subsection{Continuous Inner Products and Norms}

Consider real-valued functions $\alpha({\bf x},t)$ and $\beta({\bf x},t)$ defined over a $d$-dimensional domain, with ${\bf x}\in{\mathcal D}\subset{\mathbb R}^{d}$, and for time $t$. Particular domains of interest in the present setting are, for a given length $L$, a 1D interval ${\mathcal D}_{L}$ and the 2D square domain ${\mathcal D}_{L,L}$, defined as
\begin{equation}
\label{dldef}
{\mathcal D}_{L} = \{x\in{\mathbb R}|0\leq x\leq L\}\quad {\mathcal D}_{L,L} = {\mathcal D}_{L}\times {\mathcal D}_{L}
\end{equation}
where $\times$ indicates a Cartesian product. The $L_{2}$ spatial inner product and norm are defined as
\begin{equation}
\label{ip1ddef}
\langle \alpha,\beta\rangle_{{\mathcal D}} = \int_{{\mathcal D}}\alpha\beta d{\bf x} \quad \|\alpha\|_{{\mathcal D}} = \sqrt{\langle \alpha,\alpha\rangle_{{\mathcal D}}}\geq 0
\end{equation}
where $d{\bf x}$ is a $d$-dimensional differential element. 

For such an inner product, it is true that
\begin{equation}
\label{ipid1}
\frac{d}{dt}\frac{1}{2}\|\alpha\|^2_{{\mathcal D}} = \langle \alpha,\partial_{t}\alpha\rangle_{{\mathcal D}}
\end{equation}
where $\partial_{t}$ and $d/dt$ represent partial and total differentiation with respect to $t$, respectively.

Integration by parts (or Green's first identity) may be written as
\begin{equation}
\label{Green}
\langle \alpha, \Delta \beta\rangle_{{\mathcal D}} = -\langle \nabla\alpha, \nabla \beta\rangle_{{\mathcal D}}+\int_{\partial{\mathcal D}}\alpha {\bf n}\cdot\nabla\beta d\sigma
\end{equation}
where $\Delta$ is the Laplacian operator in $d$ dimensions, $\nabla$ is the $d$-dimensional gradient, and where $\partial{\mathcal D}$ is the boundary of the domain ${\mathcal D}$, with outward normal ${\bf n}$. 

\subsection{Discrete Domains and Grid Functions}

An $N+1$ point discrete 1D spatial domain $d_{N}$ corresponding to the line segment ${\mathcal D}_{L}$ may be defined as
\begin{equation}
d_{N} = \{l\in{\mathbb Z},0\leq l\leq N\}
\end{equation}
for some integer $N$ such that the grid spacing $h=L/N$. In a finite difference setting, other domains which differ from $d_{N}$ at the endpoints are also of interest:
\begin{equation}
d_{\underline{N}} = d_{N}-\{N\}\qquad d_{\overline{\underline{N}}} = d_{N}-\{0,N\}
\end{equation}
A 1D grid function $u_{l}^{n}$, defined over $l\in {d}$, for some 1D domain $d$, and for integer $n$ represents an approximation to $u(x,t)$ at $t=nk$ and $x=lh$. 

Similarly, a 2D discrete domain $d_{N,N}$ corresponding to the square region ${\mathcal D}_{L,L}$ may be written as
\begin{equation}
d_{N,N} = d_{N}\times d_{N}
\end{equation}
and truncated domains follow as
\begin{equation}
d_{\underline{N},N} = d_{\underline{N}}\times d_{N}\quad d_{N,\underline{N}} = d_{N}\times d_{\underline{N}}
\end{equation}
A 2D grid function $u_{l,m}^{n}$, defined over $(l,m)\in {d}$, for some 2D domain $d$, and for integer $n$ represents an approximation to $u(x,y,t)$ at $t=nk$, $x=lh$ and $y=mh$. 

\subsection{Difference Operators}

The time difference and averaging operators presented in Section \ref{timeopsec} are unchanged in their application to grid functions. 

In 1D, unit rightward and leftward spatial shifts $e_{x+}$ and $e_{x-}$ as applied to a grid function $u_{l}^{n}$ are defined as
\begin{equation}
e_{x+}u^{n}_{l} = u^{n}_{l+1}\qquad e_{x-}u^{n}_{l} = u^{n}_{l-1}
\end{equation}
Forward, backward and centered difference approximations to a first spatial derivative may be defined in terms of these shifts as
\begin{equation}
\label{dxdef}
\hspace{-0.3in}\delta_{x+} = \frac{e_{x+}-1}{h}\quad \delta_{x-} = \frac{1-e_{x-}}{h}\quad \delta_{x\cdot} = \frac{e_{x+}-e_{x-}}{2h}
\end{equation}
and centered approximation to second and fourth spatial derivatives as
\begin{subequations}
\begin{eqnarray}
\label{dxxdef}
\hspace{-0.3in}\delta_{xx} &=& \frac{e_{x+}-2+e_{x-}}{h^2}\\
\hspace{-0.3in}\delta_{xxxx} &=& \frac{e_{x+}^{2}-4e_{x+}+6-4e_{x-}+e_{x-}^2}{h^4}
\end{eqnarray}
\end{subequations}
Averaging operators may be defined as
\begin{equation}
\label{muxdef}
\hspace{-0.3in}\mu_{x+} = \frac{e_{x+}+1}{2}\quad \mu_{x-} = \frac{1+e_{x-}}{2}\quad \mu_{xx} = \mu_{x+}\mu_{x-}
\end{equation}

In 2D, unit shifts $e_{x+}$, $e_{x-}$, $e_{y+}$, and $e_{y-}$ as applied to a grid function $u_{l,m}^{n}$ are defined as
\begin{eqnarray}
e_{x+}u^{n}_{l,m} &=& u^{n}_{l+1,m}\qquad e_{x-}u^{n}_{l,m} = u^{n}_{l-1,m}\\
e_{y+}u^{n}_{l,m} &=& u^{n}_{l,m+1}\qquad e_{y-}u^{n}_{l,m} = u^{n}_{l,m-1}
\end{eqnarray}
and first difference operators in both $x$ and $y$ may be defined in analogy with \eqref{dxdef} and second difference operators $\delta_{xx}$ and $\delta_{yy}$ with \eqref{dxxdef}. A simple approximation to the Laplacian operator follows as
\begin{equation}
\label{lapfddef}
\delta_{\Delta} = \delta_{xx}+\delta_{yy}
\end{equation}
\subsection{Inner Products}

For 1D grid functions $\alpha_{l}^{n}$ and $\beta_{l}^{n}$, defined over a domain $d$, an $L_{2}$ inner product and norm may be defined as
\begin{eqnarray}
\label{ipddef}
\hspace{-0.2in}\langle \alpha,\beta\rangle_{d} &=& \sum_{l\in d}h\alpha_{l}^{n}\beta_{l}^{n} \quad \|\alpha\|_{d} = \sqrt{\langle \alpha,\alpha\rangle_{d}}\geq 0
\end{eqnarray}
For 2D grid functions $\alpha_{l,m}^{n}$ and $\beta_{l,m}^{n}$, defined over a domain $d$, the inner product and norm may be defined as
\begin{equation}
\label{ipddef}
\hspace{-0.3in}\langle \alpha,\beta\rangle_{d} = \!\!\sum_{(l,m)\in d}\!\!\!h^2\alpha_{l,m}^{n}\beta_{l,m}^{n} \quad \|\alpha\|_{d} = \sqrt{\langle \alpha,\alpha\rangle_{d}}\geq 0
\end{equation}

\subsection{Identities and Bounds}

The following identities hold for grid functions $u$ defined over a domain $d$:
\begin{eqnarray}
\langle u,\delta_{t\cdot} u\rangle_{d} &=& \delta_{t-}\frac{1}{2}\langle u,e_{t-} u\rangle_{d}\\
\langle \delta_{t\cdot}u,\delta_{tt} u\rangle_{d} &=& \delta_{t-}\frac{1}{2}\| \delta_{t+}u\|_{d}^2
\end{eqnarray}
where here, the spatial and temporal indices have been suppressed. Furthermore, the following inequality holds:
\begin{equation}
\label{dtbound}
\langle u,e_{t-} u\rangle_{d}\geq -\frac{k^2}{4}\|\delta_{t-}u\|_{d}^2
\end{equation}
Various discrete counterparts to integration by parts \eqref{Green} are available;  here are various forms of interest, written in terms of two grid functions $\alpha_{l}^{n}$ and $\beta_{l}^{n}$. In 1D, one has:
\begin{subequations}
\label{sp1ddef}
\begin{eqnarray}
\hspace{-0.3in}\langle \alpha^{n}, \delta_{xx}\beta^{n}\rangle_{d_{N}} &=& -\langle\delta_{x+}\alpha^{n},\delta_{x+}\beta^{n}\rangle_{d_{\underline{N}}}\\
&& -\alpha_{0}^{n}\delta_{x-}\beta_{0}^{n}+\alpha_{N}^{n}\delta_{x+}\beta_{N}^{n}\notag\\
\hspace{-0.3in}\langle \alpha^{n}, \delta_{xxxx}\beta^{n}\rangle_{d_{N}} &=& \langle\delta_{xx}\alpha^{n},\delta_{xx}\beta^{n}\rangle_{d_{\overline{\underline{N}}}}\\
&& -\alpha_{0}^{n}\delta_{x-}\delta_{xx}\beta_{0}^{n}+\delta_{x+}\alpha_{0}^{n}\delta_{xx}\beta_{0}^{n}\notag\\
&& +\alpha_{N}^{n}\delta_{x+}\delta_{xx}\beta_{N}^{n}-\delta_{x-}\alpha_{N}^{n}\delta_{xx}\beta_{N}^{n}\notag
\end{eqnarray}
\end{subequations}
In 1D, the action of spatial difference operators may be bounded as
\begin{equation}
\label{dxbound}
\|\delta_{xx}u\|_{d_{\overline{\underline{N}}}}\leq \frac{2}{h}\|\delta_{x+}u\|_{d_{\underline{N}}}\leq \frac{4}{h^2}\|u\|_{d_{N}}
\end{equation}
Similarly, in 2D, summation by parts in the two directions $x$ and $y$, may be written as
\begin{subequations}
\label{ip2ddef}
\begin{eqnarray}
&&\langle \alpha^{n},\delta_{xx}\beta^{n}\rangle_{d_{N,N}} = -\langle \delta_{x+}\alpha^{n},\delta_{x+}\beta^{n}\rangle_{d_{\underline{N},N}}\\
&& \qquad + \sum_{m=0}^{N}h^{2}\left(\alpha_{N,m}^{n}\delta_{x+}\beta_{N,m}^{n}-\alpha_{0,m}^{n}\delta_{x-}\beta_{0,m}^{n}\right)\notag\\
&&\langle \alpha^{n},\delta_{yy}\beta^{n}\rangle_{d_{N,N}} = -\langle \delta_{y+}\alpha^{n},\delta_{y+}\beta^{n}\rangle_{d_{N,\underline{N}}}\\
&& \qquad + \sum_{l=0}^{N}h^{2}\left(\alpha_{l,N}^{n}\delta_{y+}\beta_{l,N}^{n}-\alpha_{l,0}^{n}\delta_{y-}\beta_{l,0}^{n}\right)\notag
\end{eqnarray}
\end{subequations}
and also the bounds
\begin{equation}
\label{dxdybound}
\|\delta_{x+}u\|_{d_{\underline{N},N}}\quad, \quad\|\delta_{y+}u\|_{d_{N,\underline{N}}}\leq \frac{2}{h}\|u\|_{d_{N,N}}
\end{equation}
\section{The Hammer-Linear String Interaction}
\label{hammersec}
The nonlinear interaction of a striking hammer with a string has seen a good deal of investigation\cite{Hall87, Boutillon88, Chaigne}. 

Consider a stiff string in contact with a hammer striking from below, defined by
\begin{equation}
\label{waveeq}
\rho_{s}\partial_{tt}u = {\mathcal L}_{s}u+gf\qquad M\frac{d^2 u_{h}}{dt^2} = -f
\end{equation}
Here, $u(x,t)$ is the transverse displacement of the string in a single polarization, as a function of time $t$ and $x\in{\mathcal D}_{L}$ as defined in \eqref{dldef}, where $L$ is the string length. $\rho_{s}$ is linear density in kg/m, and $\partial_{tt}$ represents second partial differentiation with respect to time. $f$ = $f(t)$ is the force imparted to the string by a colliding hammer, of mass $M$ and at vertical height $u_{h}=u_{h}(t)$, and where $g = g(x)$ is a distribution selecting the region of impact of the hammer over the string. For a pointwise impact at a location $x=x_{0}$, one may use a Dirac distribution $g(x) = \delta(x-x_{0})$, but for a hammer of finite width, a normalized distribution with $\int_{{\mathcal D}_{L}}g dx = 1$ may be employed.

The linear operator ${\mathcal L}_{s}$ in \eqref{waveeq} is defined by
\begin{equation}
\label{lsdef}
\hspace{-0.3in}{\mathcal L}_{s} = T_{s}\partial_{xx}-E_{s}I_{s}\partial_{xxxx}-2\sigma_{s,0}\rho_{s}\partial_{t}+2\sigma_{s,1}\rho_{s}\partial_{txx}
\end{equation}
where $T_{s}$ is string tension in $N$, $E_{s}$ is Young's modulus in Pa, $I_{s}=\pi r_{s}^4/2$ is the moment of inertia of the string in m$^{4}$, where $r_{s}$ is the string radius, and $\sigma_{s,0}\geq 0$ and $\sigma_{s,1}\geq 0$ are parameters allowing for frequency dependent loss---when $\sigma_{s,0}=\sigma_{s,1}=0$, the system is lossless. Such a system is similar to that which has been used in models of lossy string vibration \cite{Ruiz, Chaigne}, but without recourse to higher time derivatives. Boundary conditions of various types may be considered, but for the present investigation clamped conditions of the type 
\begin{equation}
\label{clampeddef}
u = 0\qquad \partial_{x}u=0
\end{equation}
at the domain endpoints $x=0$ and $x=L$ are sufficient. 

The force $f$ depends on a measure of distance $\eta$ between the string and the hammer:
\begin{equation}
\label{etadef}
f = \frac{d\Phi}{d\eta}\qquad \eta = u_{h}-\langle g, u\rangle_{{\mathcal D}_{L}}
\end{equation}
The notation $\langle\cdot,\cdot\rangle_{{\mathcal D}_{L}}$ represents an $L_{2}$ inner product over ${\mathcal D}_{L}$, as defined in \eqref{ip1ddef}. Here again, $\Phi(\eta)$ is a potential function, typically modeled as a power law of the form of $\Phi_{K,\alpha}$ \cite{Chaigne}, with $K$ and $\alpha$ set from experiment; lossy models, using a force term of the form given in \eqref{lossyforcedef} are also used in studies of hysteresis of the felt in piano hammers.

\subsection{Energy Balance}
The time derivative of the total energy of the combined string/hammer system may be derived by taking an inner product of the first of \eqref{waveeq} with $\partial_{t}u$ over ${\mathcal D}_{L}$:
\begin{equation}
\rho_{s}\langle \partial_{t}u,\partial_{tt}u \rangle_{{\mathcal D}_{L}}= \langle \partial_{t}u,{\mathcal L}_{s}u\rangle_{{\mathcal D}_{L}}+\langle \partial_{t}u,g\rangle_{{\mathcal D}_{L}}f
\end{equation}
Using identities \eqref{ipid1} and integration by parts \eqref{Green}, and noting from \eqref{etadef} that, for a fixed distribution $g$, 
\begin{equation}
\langle \partial_{t}u,g\rangle_{{\mathcal D}_{L}}= \frac{d}{dt}\langle u,g\rangle_{{\mathcal D}_{L}} = \frac{d}{dt}\left(u_{h}-\eta\right)
\end{equation}
then the following energy balance results:
\begin{equation}
\frac{d{\mathcal H}}{dt} = -{\mathcal Q}_{s}+{\mathcal B}_{s}\Big|_{0}^{L}\quad{\rm ,}\quad {\mathcal H} = {\mathcal H}_{s}+{\mathcal H}_{h}
\end{equation}
Here, the string energy ${\mathcal H}_{s}$, hammer energy ${\mathcal H}_{h}$, power dissipated in the string ${\mathcal Q}_{s}$ and supplied at the boundary ${\mathcal B}_{s}$ are given by
\begin{subequations}
\label{hammeregy}
\begin{eqnarray}
\hspace{-0.3in}{\mathcal H}_{s} &=& \frac{\rho_{s}}{2}\|\partial_{t}u\|_{{\mathcal D}_{L}}^2\!+\!\frac{T_{s}}{2}\|\partial_{x}u\|_{{\mathcal D}_{L}}^2\!+\!\frac{E_{s}I_{s}}{2}\|\partial_{xx}u\|_{{\mathcal D}_{L}}^2\\
\hspace{-0.3in}{\mathcal H}_{h} &=& \frac{M}{2}\left(\frac{du_{h}}{dt}\right)^2+\Phi\\
\hspace{-0.3in}{\mathcal Q}_{s} &=& 2\sigma_{s,0}\rho_{s}\|\partial_{t}u\|_{{\mathcal D}_{L}}^2+2\sigma_{s,1}\rho_{s}\|\partial_{tx}u\|_{{\mathcal D}_{L}}^2\\
\hspace{-0.3in}{\mathcal B}_{s} &=& T_{s}\partial_{t}u\partial_{x}u-E_{s}I_{s}\partial_{t}u\partial_{xxx}u\\
\hspace{-0.3in}&& +E_{s}I_{s}\partial_{tx}u\partial_{xx}u+\sigma_{s,1}\rho_{s}\partial_{t}u\partial_{tx}u\notag
\end{eqnarray}
\end{subequations}
where ${\mathcal H}_{s},{\mathcal H}_{h},{\mathcal Q}_{s}\geq 0$. Under the clamped conditions \eqref{clampeddef}, ${\mathcal B}_{s}$ vanishes, and the system is strictly dissipative, i.e., $d{\mathcal H}/dt\leq 0$. As the individual terms in the energy balance are non-negative, it is again possible to arrive at bounds on $u_{h}$ and $u$ in terms of the total energy. 

\subsection{Finite Difference Scheme}

Employing the difference operators mentioned in the previous section, an approximation to \eqref{waveeq} over the grid $d_{N}$, of spacing $h$, is then
\begin{equation}
\label{waveeqfd}
\rho_{s}\delta_{tt}u^{n}_{l} = {\mathfrak l}_{s}u^{n}_{l}+g_{l}f^{n}\qquad M\delta_{tt}u_{h}^{n} = -f^{n}
\end{equation}
where $u_{l}^{n}$ is a grid function approximating the string displacement, and where $u_{h}^{n}$ and $f^{n}$ are time series approximating the hammer displacement and force respectively. $g_{l}$ is an approximation to the spatial distribution of the hammer. As in the continuous case, the grid function $g_{l}$ is chosen normalized such that $\sum_{l=0}^{N}hg_{l} = 1$.  

The operator ${\mathfrak l}_{s}$ is an approximation to ${\mathcal L}_{s}$:
\begin{equation}
\label{lsddef}
\hspace{-0.3in}{\mathfrak l}_{s} = T_{s}\delta_{xx}-E_{s}I_{s}\delta_{xxxx}-2\sigma_{s,0}\rho_{s}\delta_{t\cdot}+2\sigma_{s,1}\rho_{s}\delta_{t-}\delta_{xx}
\end{equation}
and the force $f$ may be written in terms of a discrete potential $\Phi^{n+\frac{1}{2}}$ as
\begin{equation}
\label{etastringdef}
f^{n} = \frac{\delta_{t-}\Phi^{n+\frac{1}{2}}}{\delta_{t\cdot}\eta^n}\qquad \eta^n = u_{h}^{n}-\langle g, u^{n}\rangle_{d_{N}}
\end{equation}
where $\Phi^{n+\frac{1}{2}} = \mu_{t+}\Phi(\eta^{n})$, and where $\langle \cdot,\cdot\rangle_{d_{N}}$ is a discrete 1D inner product, as defined in \eqref{ipddef}. 

Scheme \eqref{waveeqfd} accompanied by \eqref{etastringdef}, as in the case of the lumped collision, requires the solution of a scalar nonlinear equation of the form \eqref{Gdef}, in $r=\eta^{n+1}-\eta^{n-1}$, where
\begin{eqnarray}
\hspace{-0.3in}m &=& k^2\left(\frac{\|g\|_{d_{N}}^2}{\rho_{s}\left(1+\sigma_{s,0}k\right)}+\frac{1}{M}\right)\qquad a = \eta^{n-1}\\
\hspace{-0.3in}b &=& -2(u_{h}^{n}-u_{h}^{n-1})+\frac{k^2}{1+\sigma_{s,0}k}\langle g,\nu^{n}\rangle_{d_{N}}
\end{eqnarray}
where
\begin{equation}
\label{nudef}
\hspace{-0.3in}\nu_{l}^{n} = \!\left(\frac{2}{k}\delta_{t-}\!+\!\frac{T_{s}}{\rho_{s}}\delta_{xx}\!-\!\frac{E_{s}I_{s}}{\rho_{s}}\delta_{xxxx}\!+\!2\sigma_{s,1}\delta_{t-}\delta_{xx}\!\right)\!u^{n}
\end{equation}
which again possesses a unique solution. Notice in particular that, given values of $u$ and $u_{h}$ through time step $n$, $a$ and $b$ may be computed explicitly, and thus only a single scalar nonlinear equation must be solved in order to solve for $r$; once $r$ is known, $f^{n}$ may be calculated, and the scheme \eqref{waveeqfd} may be advanced to time step $n+1$. The scheme presented here is similar to that appearing in \cite{Chabassier}, though now accompanied by existence and uniqueness results, due to the result in \cite{chatziioannou13}, now extended to this case of the hammer/string interaction.

\subsection{Numerical Energy and Stability Condition}
An energy balance for scheme \eqref{waveeqfd} follows from an inner product with $\delta_{t\cdot}u$, and using summation by parts identities \eqref{sp1ddef}, as
\begin{equation}
\delta_{t-}{\mathfrak h}^{n+\frac{1}{2}} = -{\mathfrak q}_{s}^{n}+{\mathfrak b}_{s,0}^{n}+{\mathfrak b}_{s,N}^{n}
\end{equation}
where ${\mathfrak h}^{n+\frac{1}{2}}$, the total numerical energy is defined as ${\mathfrak h} = {\mathfrak h}_{s}+{\mathfrak h}_{h}$, and
\begin{subequations}
\label{hammerdegy}
\begin{eqnarray}
\hspace{-0.35in}{\mathfrak h}_{s}^{n+\frac{1}{2}} &=& \frac{\rho_{s}}{2}\|\delta_{t+}u^{n}\|_{d_{N}}^2+\frac{T_{s}}{2}\langle\delta_{x+}u^{n}, \delta_{x+}u^{n+1}\rangle_{d_{\underline{N}}}\\\notag
\hspace{-0.35in}&& \hspace{-0.3in}+\frac{E_{s}I_{s}}{2}\langle\delta_{xx}u^{n},\delta_{xx}u^{n+1} \rangle_{d_{\overline{\underline{N}}}}-\frac{\sigma_{s,1}k\rho_{s}}{2}\|\delta_{t+}\delta_{x+}u^{n}\|_{d_{\underline{N}}}^2\notag\\
\hspace{-0.35in}{\mathfrak h}_{h}^{n+\frac{1}{2}} &=& \frac{M}{2}\left(\delta_{t+}u_{h}\right)^2+\Phi^{n+\frac{1}{2}}\\
\hspace{-0.35in}{\mathfrak q}_{s}^{n} &=& 2\rho_{s}\left(\sigma_{s,0}\|\delta_{t\cdot}u^{n}\|_{d_{N}}^2+\sigma_{s,1}\|\delta_{t\cdot}\delta_{x+}u^{n}\|_{d_{\underline{N}}}^2\right)\\
\hspace{-0.35in}{\mathfrak b}_{s,0}^{n} &=& -T_{s}\delta_{t\cdot}u_{0}^{n}\delta_{x-}u_{0}^{n}+E_{s}I_{s}\delta_{t\cdot}u_{0}^{n}\delta_{x-}\delta_{xx}u_{0}^{n}\\
\hspace{-0.35in}&& \hspace{-0.3in}-E_{s}I_{s}\delta_{t\cdot}\delta_{x+}u_{0}^{n}\delta_{xx}u_{0}^{n}-2\sigma_{s,1}\rho_{s}\delta_{t\cdot}u_{0}^{n}\delta_{t-}\delta_{x-}u_{0}^{n}\notag\\
\hspace{-0.35in}{\mathfrak b}_{s,N}^{n} &=& T_{s}\delta_{t\cdot}u_{N}^{n}\delta_{x+}u_{N}^{n}-E_{s}I_{s}\delta_{t\cdot}u_{N}^{n}\delta_{x+}\delta_{xx}u_{N}^{n}\\
\hspace{-0.35in}&& \hspace{-0.3in}+E_{s}I_{s}\delta_{t\cdot}\delta_{x-}u_{N}^{n}\delta_{xx}u_{N}^{n}+2\sigma_{s,1}\rho_{s}\delta_{t\cdot}u_{N}^{n}\delta_{t-}\delta_{x+}u_{N}^{n}\notag
\end{eqnarray}
\end{subequations}

Here, under (for example) discrete clamped boundary conditions
\begin{equation}
u_{0}^{n} = \delta_{x+}u_{0}^{n} = 0\qquad u_{N}^{n} = \delta_{x-}u_{N}^{n} = 0
\end{equation}
then the system is numerically dissipative. For numerical stability, all that remains is to find a condition under which the energy function is non-negative---notice that ${\mathfrak h}_{s}^{n+\frac{1}{2}}$ is of indeterminate sign. To this end, using the inequalities \eqref{dtbound} and \eqref{dxbound}, it may be bounded as
\begin{equation}
\hspace{-0.35in}{\mathfrak h}_{s}^{n+\frac{1}{2}}\geq \!\left(\frac{\rho_{s}}{2}\!-\!\frac{T_{s}k^2}{h^2}\!-\!\frac{4E_{s}I_{s}k^2}{h^4}\!-\!\frac{2\sigma_{s,1}\rho_{s}k}{h^2}\right)\!\|\delta_{t+}u^{n}\|_{d_{N}}^2
\end{equation}
For a given time step $k$, this term is non-negative under the condition $h\geq h_{min}$, where
\begin{equation}
\label{stringcfl}
\hspace{-0.35in}h_{min}^2 = \frac{T_{s}k^2}{2\rho_{s}}+\frac{k\sigma_{s,1}}{2}+\frac{k}{2}\sqrt{\left(\frac{T_{s}k}{\rho_{s}}\!+\!4\sigma_{s,1}\right)^2\!\!\!+\!\!\frac{16E_{s}I_{s}}{\rho_{s}}}
\end{equation}
which serves as a stability condition for the complete scheme---it is equivalent to the condition arrived at through frequency domain (von Neumann) analysis \cite{Strikwerda} for the string alone, though now under nonlinear conditions, and reduces to the familiar Courant Friedrichs Lewy condition \cite{Courant28} in the case of an ideal lossless string.  
\subsection{Simulations}
\label{hammersimsec}
\begin{figure}[htbp]
\includegraphics[width=\columnwidth]{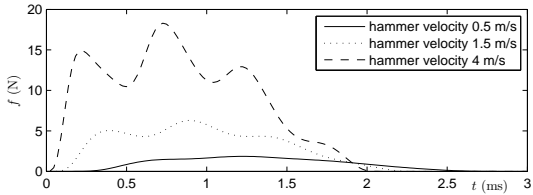}
\caption{Hammer force history $f$ in N, for a struck string, under striking velocities, as indicated. Parameters corresponding to a C4 piano string are chosen: $\rho_{s} =0.0063$ kg/m, $T_{s} = 670$ N, $E_{s}=2\times 10^{11}$ Pa, $r_{s}= 5\times 10^{-4}$ m, $\sigma_{s,0} = 0.5$ and $\sigma_{s,1}=0.5$, and $L=0.62$ m. The hammer, of mass $M= 0.0029$ kg strikes pointwise at a position 0.12 of the way along the string. The collision potential $\Phi_{K,\alpha}$ is used, with $K=4.5\times 10^{9}$ and $\alpha = 2.5$. The sample rate is 44 100 Hz.}
\label{hammerforcefig}
\end{figure}

As a simple illustration of this scheme, the case of a C4 piano string subject to a piano hammer striking action is considered here, with parameters as drawn from the article by Chaigne and Askenfelt \cite{Chaigne2}. In Figure \ref{hammerforcefig}, the force history of such a strike is plotted, under hammer velocities corresponding to piano, mezzoforte and forte strikes, illustrating the main features of a general decrease in contact time with increasing strike velocity, and the appearance of secondary humps in the force history, due to the reflection of waves from the string termination. 

\begin{figure}[htbp]

\includegraphics[width=1\columnwidth]{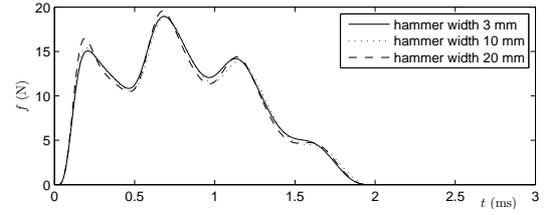}
\caption{Variations in hammer force history with hammer width. Here, the hammer interaction region $g$ is modelled as a rectangular window, of width as indicated, and the string is as in the caption to Figure \ref{hammerforcefig}, and the hammer strikes with a velocity of 4 m/s. The sample rate is 176 000 Hz.}
\label{hammerforcewidthfig}
\end{figure}

More interesting, in the present context, are more delicate features due to hammer width, which can be modeled through a distributed contact region $g$; such a choice has little impact on the computational cost of the algorithm, and subtle variations in the force history can be observed with increasing contact width, as shown in Figure \ref{hammerforcewidthfig}, where $g$ is a simple rectangular window. One issue which emerges here is of the representation of such a function on a relatively coarse grid---one solution is to operate at a high sample rate, leading to a sufficiently fine spatial grid resolution; another is to employ a high order approximation to such a finite width distribution over a coarse grid. The former approach has been adopted here.   

\begin{figure}[htbp]
\includegraphics[width=\columnwidth]{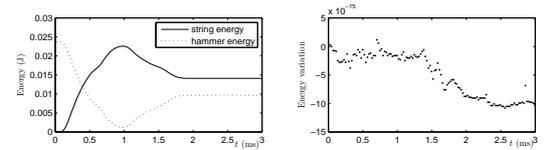}
\caption{Left: time evolution of the energy partition between the hammer and string, for the system as given in the caption to Figure \ref{hammerforcefig}, under lossless conditions (i.e., with $\sigma_{s,0}=\sigma_{s,1} = 0$). Right: normalized variation in numerical energy for the combined system.}
\label{hammerenergyfig}
\end{figure}

Numerical energy conservation is illustrated in Figure \ref{hammerenergyfig}, under lossless conditions; as in the case of the lumped collision, numerical energy is conserved to machine accuracy (here double precision floating point arithmetic). 

\begin{figure}[htbp]
\includegraphics[width=\columnwidth]{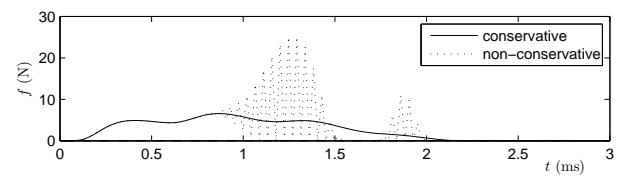}

\caption{Force history, for the system as given in the caption to Figure \ref{hammerforcefig}, under lossless and non-stiff conditions, and under a hammer velocity of 1.5 m/s, using a conservative scheme, and a nonconservative scheme. }
\label{hammerunstabfig}
\end{figure}

As a demonstration of the need for such a stability property in modeling such collisions, consider a comparison between the results of the scheme \eqref{waveeqfd} using the numerically conservative force definition from \eqref{etastringdef}, with a simpler non-conservative (and fully explicit) scheme. The simplest possible design employs a calculation of $f^{n}$ from previously computed values as
\begin{equation}
f^{n}= K[u_{h}^{n}-\langle u^{n},g\rangle_{d_{N}}]_{+}^{\alpha}
\end{equation}
Even under mild hammer excitation conditions, the non-conservative scheme exhibits severe spurious oscillations, as illustrated in Figure \ref{hammerunstabfig} under lossless and non-stiff conditions. The introduction of losses always has an ameliorating effect on such spurious oscillations, so the example described here is to be viewed as a worst case. 

\section{The Single Reed Wind Instrument}
\label{reedsec}
The single reed instrument is an interesting case of an excitation mechanism incorporating distinct nonlinearities---one is the usual nonlinear pressure/flow relation in the mouthpiece\cite{Wilson74}, but effects of collision of the reed with the lay also play a role. Reed-beating effects have been modelled by various authors---in some cases, a non-penetrative constraint is employed \cite{Guillemain04}, and in others, some penetration is allowed, and modelled through regularization \cite{karkar12, Bilbaonewbook}. A complete model requires a distributed model of the reed itself \cite{Stewart80, vanWalstijnavanzini04}; only the lumped case will be treated in this short section. A finite difference model of the reed/lay collision has been presented previously \cite{Bilbaonewbook}; in that case, however, the collision was treated using ad hoc methods, and no stability condition was available. Here, the formalism presented in the previous section, making use of a penalty potential is employed, leading to strictly dissipative numerical scheme.

Consider an acoustic tube, of length $L$, and cross-section $S(x)$, defined again over $x\in{\mathcal D}_{L}$. Here, $x=0$ corresponds to the mouthpiece termination, and $x=L$ to the bell. Lossless wave propagation in a planar mode in the tube can be described by Webster's equation\cite{Morse}:
\begin{equation}
\label{websterdef}
\frac{\rho S}{c^2}\partial_{tt}\Psi  = {\mathcal L}_{w}\Psi\qquad {\mathcal L}_{w} = \rho\partial_{x}\left(S\partial_{x}\,\cdot\,\right)
\end{equation}
where $\rho$ and $c$ are air density in kg/m$^{3}$ and sound speed in m/s, respectively, and where here, $\Psi(x,t)$ is the velocity potential in the tube. Acoustic pressure deviation $p$ and volume velocity $u$ may be derived from $\Psi$ as
\begin{equation}
p = \rho\partial_{t}\Psi\qquad u=-S\partial_{x}\Psi
\end{equation}
In particular, $p_{in} = p(0,t)$, $u_{in} = u(0,t)$, $p_{b} = p(L,t)$, $u_{b} = u(L,t)$ are the pressure/volume velocity pairs at the excitation and bell termination, respectively. In a more realistic model, viscothermal boundary layer loss effects must be included, but will be neglected here.  

The reed is modelled as a second order oscillator in displacement $z(t)$, as
\begin{equation}
\label{reedeq}
M_{r}\frac{d^2 z}{dt^2}+2M_{r}\sigma_{r}\frac{d z}{dt}+M_{r}\omega_{r}^2 z = f-S_{r}p_{\Delta}
\end{equation}
Here, $M_{r}$ is reed mass, $S_{r}$ is effective reed area, $\sigma_{r}$ a loss parameter, and $\omega_{r}$ is reed angular frequency. The force term $f$ is due to collision with the lay, located at $z=-H$, where $H$ is the equilibrium distance between the reed and lay, assumed constant here. It may be written as
\begin{equation}
\label{freeddef}
f = \frac{d\Phi}{dt}/\frac{d\eta}{dt}\qquad \eta = -z-H
\end{equation}
for some potential $\Phi(\eta)$, such as $\Phi = \Phi_{K,\alpha}$, which may be viewed as a penalty. $p_{\Delta}$ is a pressure difference across the reed, and may be written as
\begin{equation}
\label{pddef}
p_{\Delta} = p_{m}-p_{in}
\end{equation}
where $p_{m}$ is the mouth pressure supplied by the player. Furthermore, from Bernoulli's law, one has for the flow $u_{m}$ through the mouthpiece, 
\begin{equation}
\label{umdef}
u_{m} = w[-\eta]_{+}\sqrt{\frac{2|p_{\Delta}|}{\rho}}{\rm sign}(p_{\Delta})
\end{equation}
By conservation of flow, 
\begin{equation}
\label{usdef}
u_{r} = u_{m}-u_{in}
\end{equation}
where $u_{r}$, the flow induced by the reed is given by
\begin{equation}
\label{urdef}
u_{r} = S_{r'}dz/dt
\end{equation}
for an effective surface area $S_{r'}$. This may be distinct from $S_{r}$ in general, but here, for brevity, we will take $S_{r}=S'_{r}$ (if they are distinct, all the subsequent energy and stability analysis remains unchanged if \eqref{reedeq} is multiplied by the scaling factor $S_{r'}/S_{r}$---indeed, recent work indicates that these factors are indeed identical \cite{vanWalstijnavanzini07}). 

\subsection{Energy Balance}
Through an inner product of Webster's equation \eqref{websterdef} with $\partial_{t}\Psi$, and employing integration by parts \eqref{Green}, an energy balance for the acoustic tube may be written as
\begin{equation}
\label{egybalreed1}
\frac{d{\mathcal H}_{w}}{dt} = {\mathcal B}_{w}|_{0}^{L} = p_{in}u_{in}-p_{b}u_{b}
\end{equation}
where
\begin{equation}
{\mathcal H}_{w} = \frac{1}{2\rho c^2}\| \sqrt{S}\partial_{t}\Psi\|_{{\mathcal D}_{L}}^2+\frac{\rho}{2}\| \sqrt{S}\partial_{x}\Psi\|_{{\mathcal D}_{L}}^2
\end{equation}
In this simple study, effects of radiation at the bell are ignored, so $p_{b}=0$; they may easily be incorporated into the analysis presented here \cite{Bilbaonewbook}.

Similarly, by multiplying \eqref{reedeq} by $dz/dt$, and using \eqref{freeddef} to \eqref{urdef}, one may arrive at an energy balance for the reed:
\begin{equation}
\label{egybalreed2}
\frac{d{\mathcal H}_{r}}{dt}= -{\mathcal Q}_{r}-{\mathcal Q}_{m}-p_{in}u_{in}+p_{m}u_{in}
\end{equation}
where ${\mathcal H}_{r}$, ${\mathcal Q}_{r}$ and ${\mathcal Q}_{m}$, all non-negative, are given by
\begin{eqnarray}
\hspace{-0.35in}{\mathcal H}_{r} &=& \frac{M_{r}}{2}\left(\frac{dz}{dt}\right)^2+\frac{M_{r}\omega_{r}^2}{2}z^2 + \Phi\\
\hspace{-0.35in}{\mathcal Q}_{r} &=& 2M_{r}\sigma_{r}\left(\frac{dz}{dt}\right)^2\quad
{\mathcal Q}_{m} = w[-\eta]_{+}\sqrt{\frac{2}{\rho}}|p_{\Delta}|^{3/2}
\end{eqnarray}
Combining \eqref{egybalreed1} and \eqref{egybalreed2} leads to the total energy balance
\begin{equation}
\frac{d{\mathcal H}}{dt} = -{\mathcal Q}_{r}-{\mathcal Q}_{m}+p_{m}u_{in}
\end{equation}
where ${\mathcal H} = {\mathcal H}_{w}+{\mathcal H}_{r}$ is the total energy, again non-negative. Under undriven conditions, it is monotonically decreasing; under usual playing conditions (i.e., for $p_{m}\neq 0$), there is a term corresponding to power supply.

\subsection{Finite Difference Scheme}

For a grid function $\Psi_{l}^{n}$ defined over $l\in d_{N}$, a finite difference approximation to \eqref{websterdef} may be written as
\begin{equation}
\label{websterfddef}
\hspace{-0.2in}\frac{\rho \bar{S}_{l}}{c^2}\delta_{tt}\Psi_{l}^{n}  = {\mathfrak l}_{w}\Psi_{l}^{n}\qquad {\mathfrak l}_{w} = \rho\delta_{x-}\left(\mu_{x+}S\delta_{x+}\,\cdot\,\right)
\end{equation}
where $S_{l}$ is derived from the continuous bore profile $S(x)$, through sampling at locations $x=lh$, and using the averaging operation as defined in \eqref{muxdef}, and where $\bar{S}_{l} = \mu_{xx}S_{l}$. Pressure and velocity pairs at the mouthpiece and bell may be defined as
\begin{subequations}
\begin{eqnarray}
p_{in}^{n} &=& \rho\delta_{t\cdot}\Psi_{0}^{n}\qquad u_{in}^{n} = -\mu_{x-}S_{0}\delta_{x-}\Psi_{0}^{n}\\
p_{b}^{n} &=& \rho\delta_{t\cdot}\Psi_{N}^{n}\qquad u_{b}^{n} = -\mu_{x+}S_{N}\delta_{x+}\Psi_{N}^{n}
\end{eqnarray}
\end{subequations}
In the remainder of this section it is assumed that $p_{b}=0$ (so that there is no radiation), implying that $\Psi_{N}^{n} = 0$.

The reed system \eqref{reedeq} is approximated as
\begin{equation}
\label{reedfdeq}
\hspace{-0.3in}M_{r}\delta_{tt}z^{n}+2M_{r}\sigma_{r}\delta_{t\cdot}z^{n}+M_{r}\omega_{r}^2 \mu_{t\cdot}z^{n} = f^{n}-S_{r}p_{\Delta}^{n}
\end{equation}
in terms of time series $z^{n}$, $p_{\Delta}^{n}$, and $f^{n}$, defined by
\begin{equation}
\label{freedfddef}
f^{n} = \frac{\delta_{t-}\Phi^{n+\frac{1}{2}}}{\delta_{t\cdot}\eta}\qquad \Phi^{n+\frac{1}{2}} = \mu_{t+}\Phi\left(\eta^{n}\right)
\end{equation}
where $\eta^{n} = -z^{n}-H$. Equations \eqref{pddef} to \eqref{usdef} remain as written, for time series $p_{in}^{n}$, $u_{m}^{n}$ and $p_{m}^{n}$, which is assumed sampled from a given mouth pressure function $p_{m}(t)$. Equation \eqref{urdef} may be approximated as
\begin{equation}
\label{urfddef}
u_{r}^{n} = S_{r'}\delta_{t\cdot}z^{n}
\end{equation}

In contrast with the case of the hammer string interaction, due to the presence of two distinct nonlinearities, namely the Bernoulli effect and collision, there is now a pair of nonlinear equations to be solved simultaneously. To this end, note that when evaluated at $l=0$, the scheme \eqref{websterfddef} leads to an instantaneous relationship between $p_{in}^{n}$ and $u_{in}^{n}$:
\begin{equation}
p_{in}^{n} = c_{0}^{n}+c_{1}u_{in}^{n}
\end{equation}
where $c_{1}>0$ is a constant, and $c_{0}^{n}$ can be computed from values of $\Psi$ through time step $n$.  

Using this relation, in conjunction with \eqref{reedfdeq}, \eqref{freedfddef}, \eqref{urfddef}, as well as \eqref{pddef}, \eqref{umdef} and \eqref{usdef}, when viewed as relations among time series, one arrives at the pair of equations:
\begin{equation}
\label{GRsys}
G(r^{n})-gp_{\Delta}^{n}=0\qquad r^{n}-R(p_{\Delta}^{n}) = 0
\end{equation}
with $r^{n}=\eta^{n+1}-\eta^{n-1}$, for a constant $g>0$, and where
\begin{subequations}
\label{GRdef}
\begin{eqnarray}
\label{Greeddef}
\hspace{-0.35in}G(r^{n}) &=& r^{n}+b^{n}+\frac{m}{r^{n}}\left(\Phi(r^{n}+a^{n})-\Phi(a^{n})\right)\\
\label{Rreeddef}
\hspace{-0.35in}R(p_{\Delta}^{n}) &=& \!-v_{0}^{n}\!-\!v_{1}p_{\Delta}^{n}\!-\!v_{2}[-\eta^{n}]_{+}\sqrt{|p_{\Delta}^{n}|}{\rm sign}(p_{\Delta}^{n})
\end{eqnarray}
\end{subequations}
where here, the constants $m$, $v_{1}$ and $v_{2}$ are all positive, and the values $b^{n}$, $a^{n}$ and $v_{0}^{n}$ can be computed directly given known values of the state through time step $n$. As previously, under a choice of penalty potential such as $\Phi = \Phi_{K,\alpha}$, this pair of equations admits a unique solution, as shown in Appendix \ref{newtonsec}. At time step $n$, once $r^{n}$ and $p_{\Delta}^{n}$ are determined, $z^{n+1}$ and $\Psi_{0}^{n+1}$, may be calculated directly, along with the remaining values of $\Psi_{l}^{n}$. 

An energy balance follows as in the continuous case as
\begin{equation}
\delta_{t-}{\mathfrak h}^{n+\frac{1}{2}} = -{\mathfrak q}_{r}^{n}-{\mathfrak q}_{m}^{n}+p_{m}^{n}u_{in}^{n}
\end{equation}
where ${\mathfrak h}^{n+\frac{1}{2}}={\mathfrak h}_{r}^{n+\frac{1}{2}}+{\mathfrak h}_{w}^{n+\frac{1}{2}}$ and where
\begin{eqnarray}
\hspace{-0.35in}{\mathfrak h}_{w}^{n+\frac{1}{2}} &=& \frac{1}{2\rho c^2}\|\sqrt{\bar{S}}\delta_{t+}\Psi\|_{d_{N}}^2\\
\hspace{-0.35in}&& + \frac{\rho}{2}\langle\mu_{x+}S\delta_{x+}\Psi, e_{t+}\delta_{x+}\Psi\rangle_{d_{\underline{N}}}\\
\hspace{-0.35in}{\mathfrak h}_{r}^{n+\frac{1}{2}} &=& \frac{M_{r}}{2}\left(\delta_{t+}z^n\right)^2+\frac{M_{r}\omega_{r}^2}{2}\mu_{t+}\left(z^n\right)^2+\Phi^{n+\frac{1}{2}}\\
\hspace{-0.35in}{\mathfrak q}^{n}_{r} \!&=&\! 2M_{r}\sigma_{r}\left(\delta_{t\cdot}z^n\right)^2\quad{\mathfrak q}^{n}_{m} \!=\! w[-\eta^{n}]_{+}\sqrt{\frac{2}{\rho}}|p_{\Delta}^{n}|^{\frac{3}{2}} 
\end{eqnarray}
Only ${\mathfrak h}_{w}^{n+\frac{1}{2}}$ is of indeterminate sign---as in the case of the hammer/string interaction, however, it may be shown to be non-negative \cite{Bilbaonewbook} under the condition
\begin{equation}
h\geq ck
\end{equation}
which is again the Courant Friedrichs Lewy condition, now generalized to the case of Webster's equation---notice in particular that it is independent of the bore profile $S$, and is thus convenient to use in practice.

\subsection{Simulations}
 
\begin{figure}[htbp]
\includegraphics[width=1.1\columnwidth]{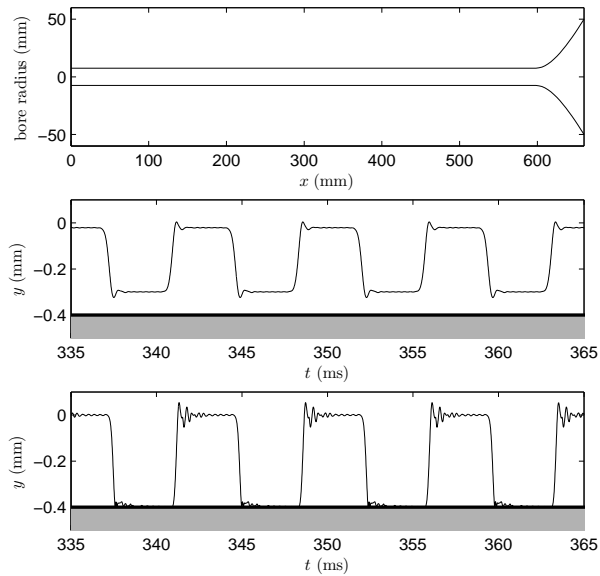}

\caption{Top: clarinet-like bore profile. Middle: steady-state oscillation of the reed position $y$, under a mouth pressure of $p_{m} = 2000$ Pa, for which the reed oscillates without beating against the lay (illustrated as a grey region at $y=-0.4$ mm). Bottom: reed position under a higher pressure of $p_{m} = 2500$ Pa, illustrating beating effects. The sample rate is 88 200 Hz.}
\label{reedposfig}
\end{figure}

\begin{figure}[htbp]
\includegraphics[width=\columnwidth]{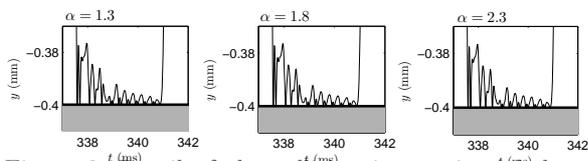}

\vspace{-0.3in}
\caption{Detail of the collision interaction of the reed under beating conditions, using a penalty potential $\Phi_{K,\alpha}$, for $K=10^{13}$ and for different values of the exponent $\alpha$, as indicated.  }
\label{reedbeatfig}
\end{figure}

As a simple test of this algorithm, consider a clarinet-like bore profile, as illustrated in the top panel of Figure \ref{reedposfig}. In this case, the reed parameters are chosen as $S_{r} = 1.46\times 10^{-4}$ m$^{2}$, $M_{r} = 3.37\times 10^{-6}$ kg, $\sigma_{r}= 1500$ s$^{-1}$, $\omega_{r} = 23250$ s$^{-1}$, and where the equilibrium distance $H$ of the reed from the lay is $4\times 10^{-4}$ m. Parameters for the acoustic field are chosen as $\rho = 1.2$ kg/m$^{3}$ and $c=340$ m/s. For the penalty potential, a choice of $\Phi = \Phi_{K,\alpha}$ is used, with $\alpha = 1.3$ and $K=10^{13}$. Both non-beating, and beating behaviour are shown; under the beating conditions, and at a sample rate of 88.2 kHz, the maximum penetration of the reed into the barrier over the duration of the simulation is 1.5$\times 10^{-8}$ m. As a check on the use of such a penalty potential for a rigid collision, plots of reed displacement under different choices of the penalty potential are shown in Figure \ref{reedbeatfig}, rather complex multiple bounce patterns are evident, the character of which are retained regardless of the choice of the exponent in the potential. 

\section{The String in Contact with a Distributed Rigid Barrier}
\label{barriersec}
In the previous examples, the colliding object, though occupying a finite region, is best considered as lumped---the nonlinear interaction force is a scalar function of time alone. In some settings, the interaction force must also be modeled as distributed. A simple example is the case of the string in contact with a rigid barrier. Though it is possible to perform a geometric analysis for the ideal string, defined by the 1D wave equation, in contact with a barrier of specified shape \cite{Burridge82, Schatzman80}, for more realistic systems, time stepping methods are probably a necessity. 

Consider a stiff and lossy string, defined again over $x\in{\mathcal D}_{L}$, and in contact with a barrier below, located at height $b(x)$. A model of this system is
\begin{equation}
\label{waveeqbarrier}
\rho_{s}\partial_{tt}u = {\mathcal L}_{s}u+{\mathcal F}
\end{equation}
where string parameters and operator ${\mathcal L}_{s}$ are as for the system described in Section \ref{hammersec}, and where ${\mathcal F}$ represents the interaction force/unit length with the barrier. If the barrier is assumed to be perfectly rigid, ${\mathcal F}$ may be written in terms of a penalty potential density $\Phi(\eta)$ as
\begin{equation}
{\mathcal F} = \partial_{t}\Phi/\partial_{t}\eta \qquad \eta = b-u 
\end{equation}
A choice of $\Phi$ of the form of a one-sided power law such as $\Phi = \Phi_{K,\alpha}$ is a convenient choice. 

Through an inner product with $\partial_{t}u$, an energy balance
\begin{equation}
\frac{d{\mathcal H}}{dt} = -{\mathcal Q}_{s}+{\mathcal B}_{s}\Big|_{0}^{L}\quad{\rm ,}\quad {\mathcal H} = {\mathcal H}{s}+{\mathcal H}_{b}
\end{equation} 
results, where ${\mathcal H}_{s}$, ${\mathcal Q}_{s}$ and ${\mathcal B}_{s}$ are as given in the case of the hammer string interaction in \eqref{hammeregy}, and where the interaction energy ${\mathcal H}_{b}$ is given by
\begin{equation}
{\mathcal H}_{b}= \int_{{\mathcal D}_{L}}\Phi dx \geq 0
\end{equation}
and the system is thus strictly dissipative. 
\subsection{Finite Difference Scheme}

A finite difference scheme, again defined over $l\in d_{N}$, follows for this system as
\begin{equation}
\label{waveeqbarrierfd}
\rho_{s}\delta_{tt}u^{n}_{l} = {\mathfrak l}_{s}u^{n}_{l}+{\mathcal F}_{l}^{n}
\end{equation}
where the difference operator ${\mathfrak l}_{s}$ is as defined in \eqref{lsddef}, and where the force density ${\mathcal F}_{l}^{n}$ is defined as
\begin{equation}
\hspace{-0.35in}{\mathcal F}_{l}^{n} \!=\! \frac{\delta_{t-}\Phi_{l}^{n+\frac{1}{2}}}{\delta_{t\cdot}\eta_{l}^{n}}\quad \Phi_{l}^{n+\frac{1}{2}}\!=\! \mu_{t+}\Phi(\eta_{l}^{n})\quad \eta_{l}^{n} \!=\! b_{l}-u_{l}^{n}
\end{equation}

The scheme satisfies an energy balance of the form
\begin{equation}
\label{barrierfd}
\delta_{t-}{\mathfrak h}^{n+\frac{1}{2}} = -{\mathfrak q}_{s}^{n}-{\mathfrak b}_{s,0}^{n}+{\mathfrak b}_{s,N}^{n}
\end{equation}
where ${\mathfrak h}^{n+\frac{1}{2}}  = {\mathfrak h}_{s}^{n+\frac{1}{2}} +{\mathfrak h}_{b}^{n+\frac{1}{2}}$ with ${\mathfrak q}_{s}^{n}$, ${\mathfrak b}_{s,0}^{n}$ and ${\mathfrak b}_{s,N}^{n}$ are as in the hammer string interaction in \eqref{hammerdegy}, and where 
\begin{equation}
{\mathfrak h}_{b}^{n+\frac{1}{2}} = \langle \Phi^{n+\frac{1}{2}},1\rangle_{d_{N}}\geq 0
\end{equation}
where ``1" indicates a grid function consisting of ones. 

Under lossless boundary conditions, the scheme is stable under the same condition as previously, namely \eqref{stringcfl}.

The scheme \eqref{barrierfd} now requires the solution of nonlinear equations along the length of the string. If $r_{l} = \eta_{l}^{n+1}-\eta_{l}^{n-1}$, then one must solve $G_{l}=0$, $l\in d_{N}$, where
\begin{equation}
G_{l} = r_{l}+\frac{m}{r_{l}}\left(\Phi(r_{l}+a_{l})-\Phi(a_{l})\right)+b_{l}
\end{equation}
with
\begin{equation}
\hspace{-0.35in}m \!=\! \frac{k^2}{\rho_{s}\left(1\!+\!\sigma_{s,0}k\right)}\quad a_{l} \!=\! e_{t-}\eta_{l}^{n}\quad b_{l} \!=\! \frac{k}{1\!+\!\sigma_{s,0}k}\nu_{l}^{n}
\end{equation}
where $\nu_{l}^{n}$ is as defined in \eqref{nudef}. In this case, existence and uniqueness follow immediately from the scalar case, as the nonlinear equations to be solved are uncoupled---which is not the case in more complex scenarios where two distributed objects are in contact. See Section \ref{snaresec}. 

In this case, where the string and barrier are assumed perfectly rigid, the penalty formulation allows some spurious penetration. This may be bounded, numerically, by noting that, because the barrier potential ${\mathfrak h}_{b}$ satisfies $0\leq {\mathfrak h}_{b}^{n+\frac{1}{2}}\leq {\mathfrak h}^{n+\frac{1}{2}}$ when the scheme is stable (i.e., under condition \eqref{stringcfl}), the penetration $\eta_{l}^{n}$ may be bounded, at all times, by
\begin{equation}
\label{stringbound}
\eta_{l}^{n}\leq \left(\frac{2\left(\alpha+1\right){\mathfrak h}^{1/2}}{Kh}\right)^{\frac{1}{\alpha+1}}
\end{equation}
Thus the penetration is bounded in terms of the initial energy, and, furthermore, can be made as small as desired through the choice of $K$. In practice, this bound is overly conservative---see the comments at the end of Section \ref{stringsimsec}.
 
\subsection{Simulations}
\label{stringsimsec}
\begin{figure}[htbp]
\includegraphics[width=\columnwidth]{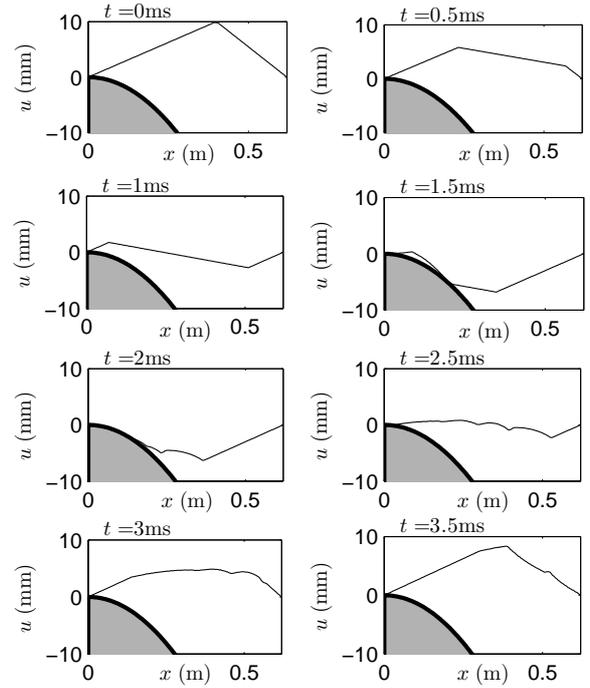}
\caption{Snapshots of the time evolution of a vibrating string in partial contact with a barrier of parabolic shape (in grey). The string is considered to be non-stiff, of length $L=0.62$ m, with $\rho_{s} = 6.3$ g/m, $T_{s} = 670$ N, $\sigma_{s,0} = 0$ and $\sigma_{1,0} = 5\times 10^{-4}$, and is initialized in a triangular shape corresponding roughly to a plucked excitation. The sample rate is 88 200 Hz. }
\label{string_evfig}
\end{figure}

The collision of a string with a rigid barrier gives rise to a wide variety of complex phenomena. As a simple example, consider a string vibrating against a rigid barrier of parabolic shape, as illustrated in Figure \ref{string_evfig}. (The parabolic shape has been used in previous studies of string collision \cite{Burridge82, kartofelev13} as a rough approximation to a bridge termination in various instruments such as the sitar or tambura; an exaggerated profile has been chosen, for illustration). In this case, the string is initialized using a triangular distribution. The potential $\Phi_{K,\alpha}$ is employed here, with $K=10^{13}$ and $\alpha = 1.3$.

Plots of the time evolution of the string are shown in Figure \ref{string_evfig}. Notice in particular the intermittent contact/recontact phenomena in evidence over the collision region, signaling that analysis or synthesis approaches based on a model of such a colliding string in terms of a moving end point may pose some difficulties \cite{Burridge82, krishnaswamy03}. 

\begin{figure}[htbp]
\includegraphics[width=\columnwidth]{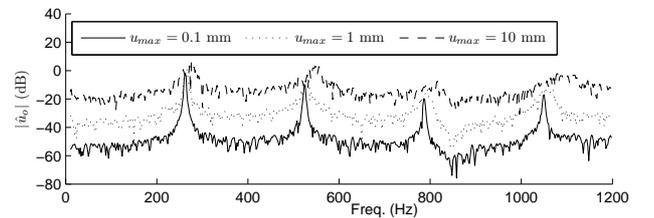}
\caption{Spectrum of string response, for the string collision as in Figure \ref{string_evfig}, under lossless conditions, and under a triangular excitation of amplitude as indicated.  }
\label{string_specfig}
\end{figure}

One perceptual effect of termination against such a smooth obstacle is an effective change in pitch with amplitude. This is illustrated in Figure \ref{string_specfig}, showing output spectra under triangular initial conditions of different amplitudes. As amplitude is increased, the spectral peaks migrate towards the higher frequencies, reflecting an effective shortening of the vibrating portion of the string, accompanied by a broadening of the peaks due to the nonlinear interaction. Under lossy conditions, one may expect pitch glide phenomena to occur. 
 
\begin{figure}[htbp]
\includegraphics[width=\columnwidth]{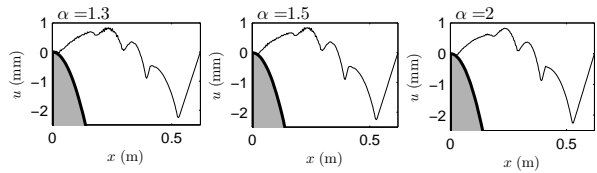}
\caption{String profile, for the string collision as in Figure \ref{string_evfig}, after 2.5 ms, under a penalty potential $\Phi_{K,\alpha}$ with $K=10^{13}$, and for different values of $\alpha$ as indicated. }
\label{stringcompfig}
\end{figure}

Interesting from a numerical perspective, given the use of a penalty based approach to rigid collision as employed here are questions of the amount of penetration, and also of convergence to an exact solution across a range of different choices of the penalty potential itself. See Figure \ref{stringcompfig}, showing the preservation of fine features of the string profile, after undergoing a collision with the barrier. For $\alpha=1.3$, and under these conditions, the maximum penetration of the string into the obstacle can be bounded, from \eqref{stringbound}, by $\eta_{l}^{n}\leq 2.6\times 10^{-5}$ m. In fact, this bound is rather conservative---for the simulation results shown in Figure \ref{string_evfig}, the maximum penetration over all grid locations, and over the length of the simulation is under 3 microns. 

\section{The Mallet-Membrane Interaction}
\label{malletsec}

The interaction of a mallet with a membrane is a generalization to 2D of the hammer string interaction described in Section \ref{hammersec}, and may be written as 
\begin{equation}
\label{malletmembrane}
\rho_{m}\partial_{tt}w = {\mathcal L}_{m}w-gf \qquad M\frac{d^2 w_{h}}{dt^2} = f
\end{equation}
Here, $w=w(x,y,t)$ is the transverse displacement of a square membrane, defined for coordinates $(x,y)\in{\mathcal D}_{L,L}$. $\rho_{m}$ is membrane density, in kg/m$^{2}$, and $g=g(x,y)$ represents the region of contact between the mallet and the membrane, again normalized such that $\int_{{\mathcal D}_{L,L}}g dxdy = 1$. $w_{h}$ represents the vertical position of the mallet, of mass $M$. The operator ${\mathcal L}_{m}$, incorporating effects of tension and simple frequency-independent loss, is defined as
\begin{equation}
\label{lmdef}
{\mathcal L}_{m} = T_{m}\Delta -2\sigma_{m,0}\rho_{m}\partial_{t}
\end{equation}
where here, $\Delta$ is the Laplacian operator, defined by
\begin{equation}
\Delta = \frac{\partial^2 }{\partial x^2}+\frac{\partial^2 }{\partial y^2}
\end{equation}
and where $T_{m}$ is the membrane tension/unit length and $\sigma_{m,0}$ is a loss parameter. The operator in \eqref{lmdef} may easily be extended, as in the case of the string, to include effects of stiffness and frequency-dependent loss; these terms will be neglected here for the sake of brevity. Boundary conditions are assumed to be of fixed type:
\begin{equation}
\label{plateclamped}
w=0 \qquad {\rm over} \qquad\partial {\mathcal D}_{L,L}
\end{equation}
where $\partial{\mathcal D}_{L,L}$ is the boundary of ${\mathcal D}_{L,L}$. 

The mallet force $f$ acting from above is defined as:
\begin{equation}
\label{eta2def}
f = \frac{d\Phi}{d\eta}\qquad \eta = \langle g, w\rangle_{{\mathcal D}_{L,L}}-w_{h}
\end{equation}
where $\langle \cdot,\cdot\rangle_{{\mathcal D}_{L,L}}$ represents a 2D inner product over the domain ${\mathcal D}_{L,L}$, and where $\Phi\geq 0$ is again a one-sided potential function, sometimes chosen\cite{Rhaouti} as a power law nonlinearity of the form $\Phi_{K,\alpha}$. 

The expression for conserved energy is directly generalized from that of the hammer string system, as:
\begin{equation}
\hspace{-0.35in}\frac{d{\mathcal H}}{dt} \!=\! -{\mathcal Q}_{m}+\int_{\partial {\mathcal D}_{L,L}}{\mathcal B}_{m}\quad{\rm ,}\quad {\mathcal H}\!=\! {\mathcal H}_{m}+{\mathcal H}_{h}
\end{equation}
Here,
\begin{subequations}
\label{memegy}
\begin{eqnarray}
{\mathcal H}_{m} &=& \frac{\rho_{m}}{2}\|\partial_{t}w\|_{{\mathcal D}_{L,L}}^2+\frac{T_{m}}{2}\|\nabla w\|_{{\mathcal D}_{L,L}}^2\\
{\mathcal H}_{h} &=& \frac{M}{2}\left(\frac{dw_{h}}{dt}\right)^2+\Phi\\
{\mathcal Q}_{m} &=& 2\sigma_{m,0}\rho_{m}\|\partial_{t}w\|_{{\mathcal D}_{L,L}}^2\\
{\mathcal B}_{m} &=& T_{m}\partial_{t}w{\bf n}\cdot \nabla w
\end{eqnarray}
\end{subequations}
where ${\bf n}$ represents a vector normal to the boundary. Under lossless boundary conditions (such as the fixed condition given in \eqref{plateclamped}), the system is dissipative. 

\subsection{Finite Difference Scheme}

The scheme, which is a direct extension of that of the hammer/string system, can be presented briefly here:
\begin{equation}
\label{waveeq2dfd}
\hspace{-0.3in}\rho_{m}\delta_{tt}w^{n}_{l,m} = {\mathfrak l}_{m}w^{n}_{l,m}-g_{l,m}f^{n}\quad M\delta_{tt}w_{h}^{n} = f^{n}
\end{equation}
where $w_{l,m}^{n}$, $w_{h}^{n}$ and $f^{n}$ are approximations to $w$ (now defined over the 2D discrete domain $d_{N,N}$), $w_{h}$ and $f$, respectively, and where ${\mathfrak l}_{w}$ is an approximation to ${\mathcal L}_{w}$:
\begin{equation}
\label{lmfddef}
{\mathfrak l}_{w} = T_{w}\delta_{\Delta}-2\sigma_{m,0}\rho_{m}\delta_{t\cdot}
\end{equation}
and the force $f$ may be written in terms of a discrete potential $\Phi^{n+\frac{1}{2}}$ as
\begin{equation}
\label{etamemdef}
\hspace{-0.3in}f^{n} = \frac{\delta_{t-}\Phi^{n+\frac{1}{2}}}{\delta_{t\cdot}\eta^n}\qquad \eta^n = \langle g, w^{n}\rangle_{d_{N,N}}-w_{h}^{n}
\end{equation}
where $\Phi^{n+\frac{1}{2}} = \mu_{t+}\Phi(\eta^{n})$. 

As previously, updating scheme \eqref{waveeq2dfd} requires the solution of a  nonlinear equation of the form \eqref{Gdef}, where
\begin{eqnarray}
\hspace{-0.3in}m &=& k^2\left(\frac{\|g\|_{d_{N,N}}^2}{\rho_{m}\left(1+\sigma_{m,0}k\right)}+\frac{1}{M}\right)\quad a = e_{t-}\eta^{n}\\
\hspace{-0.3in}b &=& 2k\delta_{t-}w_{h}^{n}\!-\!\frac{2k}{1\!+\!\sigma_{m,0}k}\langle g,\left(\delta_{t-}\!\!+\frac{kT_{m}}{2\rho_{m}}\delta_{\Delta}\right)\!w^{n}\rangle_{d_{N,N}}\notag
\end{eqnarray}

The energy balance is now, after taking an inner product over $d_{N,N}$ with $\delta_{t\cdot}w$, and using \eqref{ip2ddef},
\begin{equation}
\delta_{t-}{\mathfrak h}^{n+\frac{1}{2}} = -{\mathfrak q}_{m}^{n}+\sum_{(\zeta_{x},\zeta_{y})\in\partial d_{N,N}}{\mathfrak b}_{m, \zeta_{x},\zeta_{y}}^{n}
\end{equation}
where $\partial d_{N,N}$ indicates the set of grid points lying on the boundary of $d_{N,N}$. ${\mathfrak h}^{n+\frac{1}{2}}$, the total numerical energy is defined as ${\mathfrak h} = {\mathfrak h}_{m}+{\mathfrak h}_{h}$, where
\begin{subequations}
\label{memdegy}
\begin{eqnarray}
\hspace{-0.25in}&& \hspace{-0.25in}{\mathfrak h}_{m}^{n+\frac{1}{2}} = \frac{\rho_{m}}{2}\|\delta_{t+}w^{n}\|_{d_{N,N}}^2+\\
\hspace{-0.25in}&& \hspace{-0.25in}\frac{T_{m}}{2}\left(\langle\delta_{x+}w^{n}, \delta_{x+}w^{n+1}\rangle_{d_{\underline{N},N}}+\langle\delta_{y+}w^{n}, \delta_{y+}w^{n+1}\rangle_{d_{N,\underline{N}}}\right)\notag\\
\hspace{-0.25in}&&\hspace{-0.25in}{\mathfrak h}_{h}^{n+\frac{1}{2}} = \frac{M}{2}\left(\delta_{t+}w_{h}\right)^2+\Phi^{n+\frac{1}{2}}\\
\hspace{-0.25in}&&\hspace{-0.25in}{\mathfrak q}_{m}^{n} = 2\sigma_{m,0}\rho_{m}\|\delta_{t\cdot}w^{n}\|_{d_{N,N}}^2\\
\hspace{-0.25in}&&\hspace{-0.25in}{\mathfrak b}_{m, \zeta_{x},\zeta_{y}}^{n}=  T_{m}\delta_{t\cdot}w_{\zeta_{x},\zeta_{y}}^{n}\delta_{b+}w_{\zeta_{x},\zeta_{y}}^{n}        
\end{eqnarray}
\end{subequations}
Here, the notation $\delta_{b+}$ indicates a spatial first difference operation in the direction normal to the boundary. At corner points in the domain $d_{N,N}$, it is to be applied in both directions and summed. Under the condition that $w_{l,m}^{n}$ is zero at the boundary, corresponding to a fixed termination \eqref{plateclamped}, the system is again dissipative.

A stability condition for scheme \eqref{waveeq2dfd} follows again from the non-negativity of ${\mathfrak h}_{m}$. Using bounds \eqref{dtbound} and \eqref{dxdybound} leads to the condition
\begin{equation}
\label{cfl2d}
h\geq h_{min} =  k\sqrt{\frac{2T_{m}}{\rho_{m}}}
\end{equation}
which is the same as the bound obtained using frequency domain techniques for the membrane scheme in isolation. 

\subsection{Simulations}
\begin{figure}[htbp]
\includegraphics[width=1.1\columnwidth]{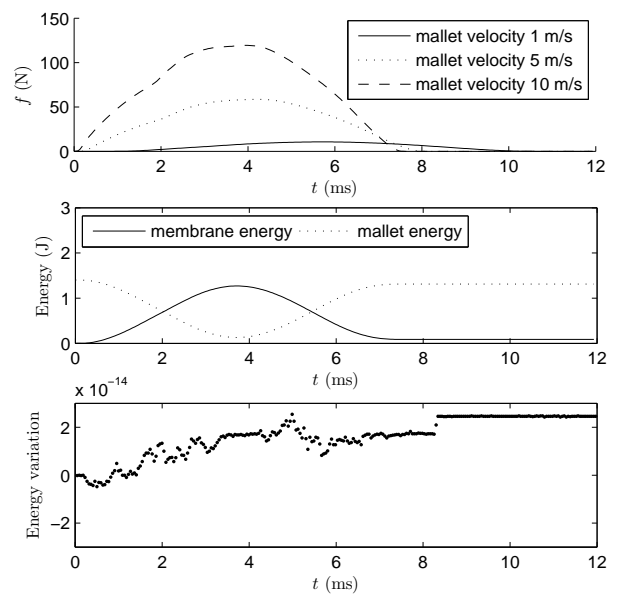}
\caption{Top: Force experienced by a mallet, of mass 0.028 kg, and with stiffness parameters $K=1.6\times 10^{8}$ and $\alpha = 2.54$, striking a membrane under different velocities, as indicated. The membrane, of dimensions 0.6$\times$ 0.6 m, with $\rho_{m} = 0.26$ kg/m$^{2}$ and $T_{m} =$ 3325 N/m, is struck at a location 0.1 m from a corner. Middle: Energy partition. Bottom: normalized energy variation. The sample rate is 22 050 Hz.}
\label{malletforcefig}
\end{figure}

In this section, the results of a simulation for a square membrane are shown, corresponding roughly to a typical drum configuration. Force histories are shown in Figure \ref{malletforcefig} at top, illustrating again the decrease in contact duration with mallet velocity, and at middle and bottom, the energy partition between the mallet and membrane, as well as the normalized energy variation are shown. The variation in the energy is compounded by the number of degrees of freedom of the system, which is considerably larger than in the previous cases. 

\section{Fully Distributed Interaction}
\label{snaresec}

The methods described in the previous sections can also be applied  
to more complex scenarios involving distributed/distributed coupling.  
In this section, the collision of an ideal string with a membrane is  
described, with an eye towards full simulation of the snare drum,  
which will be briefly outlined.

\subsection{The String/Membrane Collision}
Consider the system of an interacting string and membrane, defined over the regions  
$\lambda \in \mathcal{D}_{L_s}$ and $\mathbf{x} = (x,y)\in \mathcal{D}_{L_m,L_m}$, and with transverse displacements $u=u(\lambda,t)$ and $w=w(\mathbf{x},t)$. The equations of motion are:
\begin{equation}\label{eq:str_memb}
\hspace{-0.35in}\rho_m \D_{tt} w = \mathcal{L}_m w + \left< g, \F  
\right>_{\mathcal{D}_{L_s}} \,\,\,\, %
\rho_s \D_{tt} u = \mathcal{L}_s u - \F
\end{equation}
The operators ${\mathcal L}_{m}$ and ${\mathcal L}_{s}$  
have been introduced in \eqref{lmdef} and \eqref{lsddef}, respectively. The  
collision force $\F$ is again a distributed function  
$\F~=~\F(\lambda)$, and the contribution of $\F$ to the  
membrane equation must be integrated along the domain of interaction with  
the string with a suitable distribution function. The simplest  
possible choice is probably $g=\delta(\mathbf{x}-\pi(\lambda))$, where  
$\pi: \mathcal{D}_{L_s} \rightarrow \mathcal{D}_{L_m,L_m}$ is the  
projection from each point of the string to the corresponding point on  
the membrane and $\delta$ is a 2D Dirac delta function.
As before, $\F$ may be related to a potential  
$\Phi=\Phi(\eta(\lambda))$, which depends on the distance $\eta$  
between the string and membrane over the region of interaction:
\begin{equation}\label{eq:fcoll_memb}
\F=\frac{\D_t \Phi}{\D_t \eta}, \qquad \eta = u- \left<g,  
w\right>_{\mathcal{D}_{L_m,L_m}}
\end{equation}

Once again, one can derive an expression for energy which is conserved  
in the lossless case.

\subsection{Finite Difference Scheme}
A finite difference approximation for (\ref{eq:str_memb}) and (\ref{eq:fcoll_memb})
can be written as follows:
\begin{subequations}
\label{eq:fd_str_memb}
\begin{align}
\rho_m \DD_{tt} w_{l,m}^n &= \frakl_m w_{l,m}^n + ({\mathfrak i}_{sm} \F^n)_{l,m}\\
\rho_s \DD_{tt} u_{l}^n &= \frakl_s u_{l}^n - \F_l^n
\end{align}
\end{subequations}
in terms of grid functions $u_{l}^{n}$, defined over $d_{N_s}$, and $w_{l,m}^{n}$ over $d_{N_m, N_m}$, where $\F_{l}^{n}$, may be written as
\begin{equation}
\label{eq:fd_f}
\F_l^n = \frac{\DD_{t-} \Phi_l^{n+\frac{1}{2}}}{\DD_{t\cdot} \eta_l^n}, \quad \eta_l^n =  
u_l^n - ({\mathfrak i}_{ms}w^n)_l
\end{equation}
Here, integrals have been expressed as linear operators:
\begin{equation}
\left<g, \,\cdot \, \right>_{d_{N_s}}= {\mathfrak i}_{sm}, \qquad \left< g,  
\,\cdot \,\right>_{d_{N_m,N_m}} ={\mathfrak i}_{ms}
\end{equation}

Ultimately, the finite difference scheme can be updated by solving a  
non-linear equation in the vector ${\bf r}$:
\begin{equation}
\label{Gvecdef}
{\bf G}({\bf r}) ={\bf r}+{\bf M}\,{\bm \omega}+{\bf b} = {\bf 0}
\end{equation}
with the elements of ${\bf r}$ and ${\bm\omega}$ given by
\begin{equation}
\hspace{-0.3in}r_l=\eta_l^{n+1}-\eta_l^{n-1}, \quad \omega_l=\frac{\Phi(r_l+a_l)-\Phi(a_l)}{r_l}.
\end{equation}
As before, $a_l=e_{t-} \eta^n_l$ and ${\bf b}$ depends only on known values of $w$ and $u$. In contrast with the previous cases, values of the solution ${\bf r}$ are now coupled by the presence of a square matrix ${\bf M}$ defined as:
\begin{equation}
\label{Mdef}
\hspace{-0.3in}{\bf M} = \left(\frac{k^2}{\rho_s(1\!+\!\sigma_{s,0}k)}\, \mathbf{1} +  
\frac{k^2}{\rho_m (1\!+\!\sigma_{m,0}k)}\, \I_{ms} \I_{sm}\right)
\end{equation}
where $\mathbf{1}$ represents the identity matrix.  $\I_{sm}$ and $\I_{ms}$ are the matrix forms of the operators ${\mathfrak i}_{sm}$ and ${\mathfrak i}_{ms}$, respectively, and are required to be the transposes of one another for energy conservation reasons\cite{Bilbaonewbook}. Therefore, ${\bf M}$ is positive definite, which guarantees existence and uniqueness of the solution (see Appendix \ref{newtonsec}.)

\subsection{Boundary conditions}
Boundary conditions for the membrane have been discussed in Section  
\ref{malletsec}. For the string, one possibility is to use fixed  
termination at both ends; a more realistic choice, when dealing with musical  
instruments, is to attach the string directly to the membrane surface.  
This design resembles that of a snare drum, where two bridges connect  
both ends of a set of snares to the membrane, while remaining in  
direct contact with it.
In this case, additional force terms resolved at the boundaries of the  
string must be added to the equation for the membrane:
\begin{eqnarray}
\hspace{-0.2in}\rho_m \DD_{tt} w_{l,m}^n &=& \frakl_m w_{l,m}^n + g^{(0)}_{l,m} f^{(0)} + g^{(L_{s})}
_{l,m} f^{(L_s)}\\
\hspace{-0.2in}&&  +  \left({\mathfrak i}_{sm} \F^n\right)_{l,m}\notag
\end{eqnarray}
where $g^{(0)}$ and $g^{(L_s)}$ are suitable distribution functions (e.g.,  
2D Dirac delta functions) for $\lambda=0, {L_s}$.
Energy analysis can be applied to find stable boundary conditions at  
$\lambda=0$ (and similarly at $\lambda={L_s}$):
\begin{multline}
\left< g^{(0)}, w \right>_{d_{N_m,N_m}} = u_0 \\
 f^{(0)} = (T_s \, \DD_{\lambda-}  - E_s  
I_s\DD_{\lambda-}\DD_{\lambda\lambda} +  
2\sigma_{s,1}\rho_{s}\DD_{t-}\DD_{\lambda-}) u_0
\end{multline}
The first condition is that displacements are equal at the snare endpoint and the point on the membrane to which is connected, and the second that the force exerted on the membrane is equal and opposite to that exerted on the end of the snare. Under these conditions, the scheme as a whole is dissipative, and stable under the separate conditions \eqref{stringcfl} and \eqref{cfl2d}. The expression for matrix ${\bf M}$ in \eqref{Mdef} is slightly altered, but remains positive definite.   

\subsection{Case Study: The Snare Drum}
In this section, an application to the simulation of the snare drum is briefly outlined. The snare drum is double-headed, with the distinguishing feature of a set of snares (metal wires) in contact with the lower membrane, which gives this  
instrument its characteristic rattling sound.

Finite difference time domain simulation of the snare drum has been presented previously\cite{Bilbaosnare}; in that case, however, the collision of the snares with the membranes was treated using ad hoc methods, and no stability condition was available. Here, again, the penalty potential formalism described in the previous sections is extended to handle this more complex case, leading to a provably stable algorithm.

\begin{figure}[htbc]
\vspace{-0.2in}
\includegraphics[width=1\columnwidth]{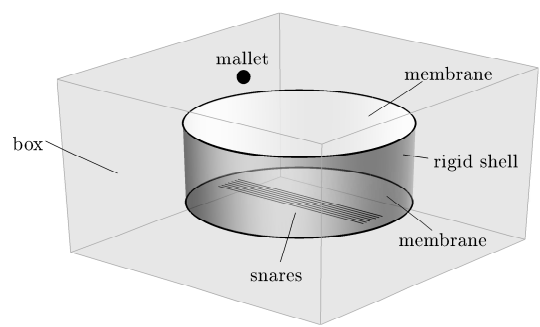}
\vspace{-0.1in}

\caption{Diagram of snare drum geometry, illustrating the various  
components. A finite enclosure with absorbing conditions at the walls  
surrounds the system.}
\label{fig:sdscheme}
\end{figure}

A simple model of the snare drum consists of two circular membranes with fixed boundary conditions, coupled by a rigid cavity and immersed in a box of air, with a set of snares in  
contact with the lower membrane (see Figure \ref{fig:sdscheme}). The  
system can be excited through a mallet acting on the upper membrane.  
Absorbing conditions are applied at the walls of the enclosure, ideally simulating an anechoic space.

\begin{figure}[htbc]
\includegraphics[width=1\columnwidth,trim=0.5cm 0cm 0.5cm 0cm,clip]{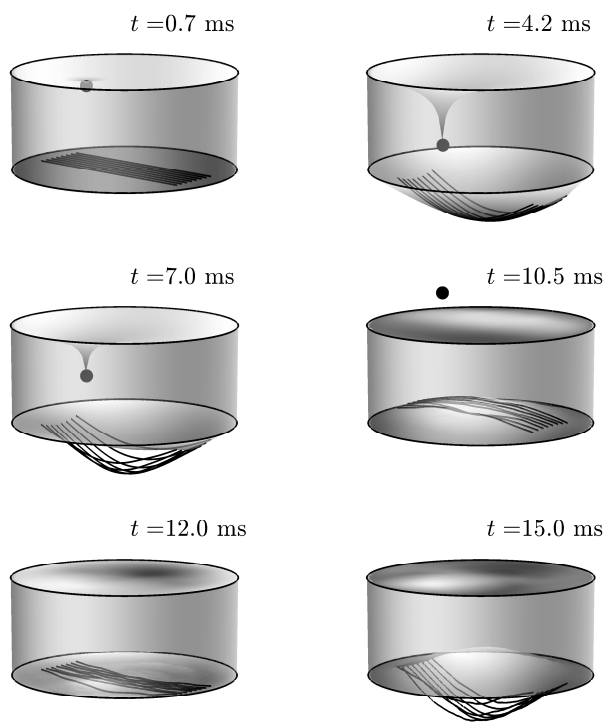}

\caption{Snapshots of the evolution of the snare drum system at times  
as indicated. Displacements have been scaled for illustration  
purposes. The sample rate is 48 000 Hz.}
\label{fig:sdsnaps}
\end{figure}

\begin{figure}[htbc]
\vspace{-0.2in}
\includegraphics[width=1\columnwidth,trim=0 0 1cm 0,clip]{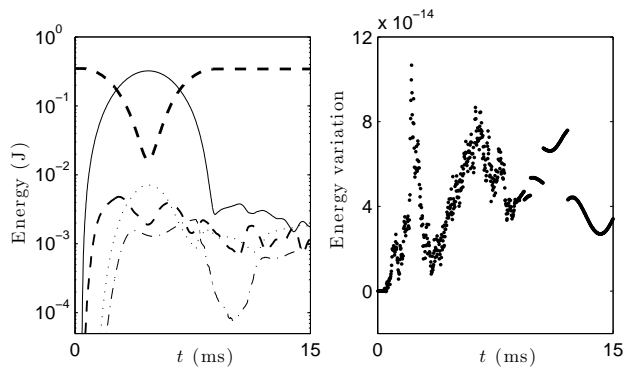}
\vspace{-0.1in}
\caption{Left: energy partition among the various components of the  
snare drum (solid line: upper membrane, bold dashed: mallet, dotted:  
lower membrane, dashed: air, dot-dashed: snares). Right: normalized  
energy variations.}
\label{fig:sdenergy}
\end{figure}

The approach described above may be applied to the  
mallet-membrane interaction and to the snares colliding against the  
membrane, in order to obtain a fully energy conserving model. Figure  
\ref{fig:sdsnaps} shows some snapshots of the time evolution of the system, subject to an  
initial strike on the upper membrane. The delay in the excitation caused by the air inside the cavity on the lower membrane is apparent. The snares are launched by the contact with the lower membrane; though the movement is at first coherent, it is rapidly randomized through multiple collisions. Figure \ref{fig:sdenergy} shows the contributions to total energy of the various components, with normalized energy variations on the order of machine accuracy.

\section{Concluding Remarks}

This article has explored basic features of time domain numerical simulation of collision interactions in musical instruments, with a focus on a Hamiltonian formulation, within which the collision mechanism may be included through an added potential; such a formulation allows for the stable simulation of a wide range of collision interactions in musical instruments. It is hoped that such techniques will aid in the investigation of finer features in musical instrument acoustics. From a sound synthesis perspective, one issue which has not been touched upon here is aliasing---always present in nonlinear models, but somewhat alleviated by physical damping. 

When the colliding objects are assumed to be rigid, the potential has the interpretation of a penalty. As has been shown at various instances, a penalty formulation of rigid collisions does lead to spurious penetration, but generally this is very small by the standard of acoustics applications---furthermore, there is some degree of control over the amount of penetration, which may be conveniently bounded, provided a numerical energetic framework is available. One aspect worthy of further study is the particular choice of such a penalty; a power law has been used here, for simplicity, but many others are available, perhaps with superior properties in terms of the reduction of spurious penetration.

The formulations presented here require the solution of nonlinear equations in the main update. If the colliding object is lumped, as in the case of the hammer/string, mallet/membrane and reed interactions then a single equation to be solved results. For fully distributed collisions between two deformable objects, the solution to a system of coupled equations is required. In all cases examined here, however, the system to be solved possesses a unique solution---see Appendix \ref{newtonsec}. The Newton Raphson method has been employed here in order to solve such nonlinear equations. One aspect of such iterative methods which has not been addressed here is convergence, even though a unique solution exists. The Newton Raphson method employed here has in all cases led to convergence (easily observed by energy conservation to machine accuracy), but as yet its convergence remains unproven for the systems examined here. Some partial results, however, are available---see Appendix \ref{newtonsec}. Newton Raphson, however, is but one method of solution---many others are available \cite{Kelley}. Useful also, here, would be bounds on the number of required iterations. 

Finally, though the Hamiltonian framework presented here has been applied to the case of finite difference time domain methods, one might suspect that it applies more generally in different formulations which ultimately reduce to time stepping methods (including time domain finite element methods, spectral methods, as well as modal techniques), though over more general basis functions. Various references are available\cite{Evans, Quarteroni} giving an overview of such families of numerical techniques. Such choices will form the basis of future investigations.

\appendix


\section{Nonlinear Equations: Existence/Uniqueness and Convergence}
\label{newtonsec}

The solution of the nonlinear equation $G(r)=0$, for $G(r)$ as defined in \eqref{Gdef}, as well as several variants, plays a key role in all of the algorithms described here. It is thus worth outlining here some properties of such equations, especially with regard to the existence and uniqueness of solutions, as well as the convergence of the Newton Raphson method. 

It has recently been shown \cite{chatziioannou13} that the equation $G(r)=0$, for $G(r)$ as defined in \eqref{Gdef} has a unique solution in $r$. The essence of the proof is that for a differentiable and convex potential $\Phi$, $G'(r)$ is non-negative, and bounded away from zero, and thus $G(r)=0$ possesses a unique solution. In particular, $G'(r)\geq 1$. The power law potential $\Phi_{K,\alpha}$ is convex and differentiable, and thus the numerical methods here that rely on the solution of such an equation, namely the collisions described in Sections \ref{lumpedsec}, \ref{hammersec}, \ref{barriersec} and \ref{malletsec}, admit a unique update. 

The Newton Raphson method can be shown to be globally convergent, independent of the starting point $r^*$, if $G(r)$ is itself convex \cite{Ortega}. Indeed, $\forall \alpha\geq2$, $G''(r) = \frac{m}{r^3}P(r)$ where $P(r) = r^2\Phi''(r+a)-2r\Phi'(r+a)+2\Phi(r+a)-2\Phi(a)$. Now, $P'(r)=r^2\Phi'''(r+a)\geq0$, since $\Phi_{K,\alpha}$ is convex, hence $P(r)$ is a monotonically increasing function passing through the origin, which results in $G''(r)\geq0$ (which also holds in the limit $r\to0$). When $1<\alpha<2$, $\Phi'''$ ceases to exist for certain values of $r$ and $a$. In this case only local convergence can be guaranteed, subject to a good initial guess $r^*$, which is ensured in practice by choosing the previous value of $r$. 

The case of the beating reed requires the solution of a pair of nonlinear equations, as given in \eqref{GRsys}, with the functions $G$ and $R$ as defined in \eqref{GRdef}. First, note that $R$ is one-to-one, thus possessing an inverse $R^{-1}$, so system \eqref{GRsys} may be condensed to 
\begin{equation}
J(r^{n}) = G(r^{n})-gR^{-1}(r^{n}) = 0
\end{equation}
Note that $R'$ is negative, and bounded between $-v_{1}$ and $-\infty$. and thus $(R^{-1})'$ is also negative, and bounded between 0 and $-1/v_{1}$. Because $g>0$, and $G'\geq 1$, $J'\geq 1$, and the system possesses a unique solution. It is not, however, convex, and a simple conclusion regarding the convergence of Newton Raphson is not  easily arrived at.

Finally, consider the vector system of nonlinear equations arising in the case of distributed/distributed interaction in \eqref{Gvecdef}. Because ${\bf M}>0$, the system to be solved is equivalent to
\begin{equation}
\hat{{\bf G}}({\bf r}) = {\bf M}^{-1}{\bf r}+{\bm \omega}+{\bf M}^{-1}{\bf b} = {\bf 0}
\end{equation}
The system possesses a unique solution if the Jacobian of $\hat{{\bf G}}$ is positive definite. The Jacobian is ${\bf M}^{-1}+{\bm \Omega}$, where ${\bm\Omega}$ is a diagonal matrix, the diagonal entries of which are the derivatives of the components of ${\bm \omega}$ which are individually non-negative, and is thus positive semi-definite. Thus the Jacobian is positive definite.

\subsection*{Acknowledgement}

This work was supported by the European Research Council, under grant number StG-2011-279068-NESS. 

\References{}{master}{actalit}

\end{document}